\documentclass[final]{siamltex}
\usepackage{color}

\usepackage{graphicx,amssymb,amsmath}
\usepackage{empheq}
\usepackage{pifont}% http://ctan.org/pkg/pifont
\newcommand{\cmark}{\ding{51}}%
\usepackage[tableposition=top]{caption}

\newcommand{\kchg}[1]{{\color{blue} #1}}

\newcommand{\mc}{\mathcal}
\newcommand{\mb}{\mathbf}

\def\e1{{\varepsilon_{11}}}
\def\b1{{\beta_{11}}}
\def\bp3{{\beta_{33}}}
\def\ep3{{\varepsilon_{33}}}
\def\Re{{\rm Re \, }}
\def\mX{{\mathbb X}}
\def\Ltwo{{\mathbb L}^2 }

\MHInternalSyntaxOn
\providecommand*\phantomword[3][c]{%
\mathchoice
{\MT_phantom_word:NNnn #1\displaystyle {#2}{#3}}%
{\MT_phantom_word:NNnn #1\textstyle {#2}{#3}}%
{\MT_phantom_word:NNnn #1\scriptstyle {#2}{#3}}%
{\MT_phantom_word:NNnn #1\scriptscriptstyle {#2}{#3}}%
}
\def\MT_phantom_word:NNnn #1#2#3#4{%
\@begin@tempboxa\hbox{$\m@th#2#4$}%
\setlength\@tempdima{\widthof{$\m@th#2#3$}}%
\hbox{\hb@xt@\@tempdima{\csname bm@#1\endcsname}}%
\@end@tempboxa}
\MHInternalSyntaxOff

\title{Modeling and stabilizability of voltage-actuated piezoelectric beams with magnetic effects}
\author{K. A. Morris  \thanks{Department of Applied Mathematics, University of Waterloo, Waterloo, ON N2L3G1, Canada ({\tt kmorris@uwaterloo.ca}).}
        \and A. \"{O}. \"{O}zer \thanks{({\tt aozer@uwaterloo.ca}).}}
\begin{document}

\maketitle

\begin{abstract}
Models for piezoelectric beams and structures with piezoelectric patches generally ignore magnetic effects.
This is because the magnetic energy  has a relatively small effect on the overall dynamics.  Piezoelectric beam models are known to be exactly observable, and can be exponentially stabilized in the energy space by using a mechanical feedback controller. In this paper, a variational approach is used  to derive a model for a piezoelectric beam that includes magnetic effects.  It is proven that the partial differential equation model is well-posed. Magnetic effects have a strong effect on the stabilizability of the control system.    For almost all system parameters the piezoelectric beam  can be strongly stabilized, but is not exponentially stabilizable in the energy space. Strong stabilization is achieved using only electrical feedback.  Furthermore, using  the same electrical feedback,  an exponentially stable closed-loop system can be obtained for a  set of system parameters of zero Lebesgue measure. These results are compared to those of a beam without magnetic effects.
 \end{abstract}


%\begin{AMS}
%15A15, 15A09, 15A23
%\end{AMS}

\pagestyle{myheadings}
\thispagestyle{plain}
\markboth{
%DRAFT-SIAM Journal of Control and Optimization
}{%OZER AND MORRIS
}

 \enlargethispage*{3pc}

  \section{Introduction}
  Piezoelectric actuators have a unique characteristic of converting mechanical energy to electrical and \emph{magnetic energy}, and vice versa. Therefore they could be used as actuators or sensors.
  Piezoelectric actuators are generally scalable, smaller, less expensive and more efficient than traditional actuators, and hence, a competitive choice for many tasks in  industry, particularly those involving control of structures.
 Piezoelectric materials been employed in civil, industrial, automotive, aeronautic, and space structures.

 %Many different PDE (Partial Differential Equations) models which incorporate the physics are developed in the context of infinite dimensional parameter estimation
%and control formulations.

%In classical mechanics, it is very well-known that equations of motion can be formulated either
%through a set of differential equations, or through a variational principle, so-called
%Hamilton's principle. In applying Hamilton's principle, the functional is specified over a fixed time interval,
%and the admissible variations of the generalized coordinates (independent variables) are taken to be
%zero. The set of field equations for the piezoelectric beams/plates have been well-established through
%the combination of beams/plates equations and Maxwell's equations.
\begin{figure}[h!tb]
\centering
\includegraphics[width=2.9in]{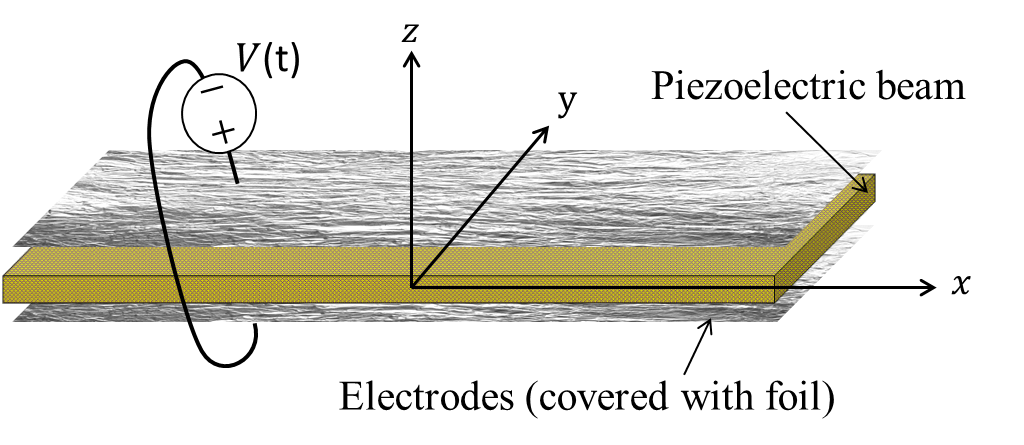}
\caption{{\footnotesize{For a voltage-actuated beam/plate, when voltage $V(t)$ is supplied to the electrodes, an electric field is created between the electrodes, and therefore the beam/plate either shrinks or extends.}}}
\label{pbeam}
\end{figure}

In modeling  of piezoelectric systems, three major effects and their interrelations need to be considered: mechanical, electrical, and magnetic.  Mechanical effects are generally modeled through Kirchhoff, Euler-Bernoulli, or Mindlin-Timoshenko small displacement assumptions; see, for instance, \cite{Banks-Smith}, \cite{Hansen}, \cite{Smith}, \cite{Yang}. To include electrical and magnetic effects, there are mainly three approaches: \emph{electrostatic, quasi-static,} and \emph{fully dynamic} \cite{Tiersten}.  Electrostatic and quasi-static approaches are widely used  - see, for instance, \cite{Dest}, \cite{Hansen}, \cite{K-M-M}, \cite{L-M}, \cite{Rogacheva}, \cite{Smith}, \cite{Tiersten}, \cite{Tzou}. These models completely exclude magnetic effects and their coupling  with electrical and mechanical effects. In a electrostatic approach,  electrical effects are  stationary, even though  the mechanical equations are dynamic.  %In other words, electrostatic theory assumes that electrical effects moves from one point to another instantaneously.
%  This approach is not far from reality since it is observed that magnetic effects do not have much to do with the dynamics of these structures.
In the case of quasi-static approach, magnetic effects are still ignored but  electric charges  have time dependence. The electromechanical coupling is  not dynamic.

A piezoelectric beam  is an elastic beam  with  electrodes at its top and bottom surfaces, insulated at the edges (to prevent fringing effects), and connected to an external electric circuit. (See Figure \ref{pbeam}). These are the simplest structures on which to study the interaction between the electrical and mechanical  energy in these systems.
It is experimentally observed that  the magnetic effects are minor in the overall dynamics for polarized ceramics  (see the review article \cite{Yang1}), and therefore these effects are  ignored in piezoelectric beam models.
 A single piezoelectric beam either shrinks or extends when the electrodes are subjected to a voltage source. For a beam of length $L$ and thickness $h,$  models derived by electrostatic and quasi-static approaches with the Euler-Bernoulli small displacement assumptions (no damping)  describe the stretching motion as
\begin{subequations}
  \label{or}
\begin{empheq}[left={\phantomword[r]{0}{ }  \empheqlbrace}]{align}
  & \rho v_{tt}-\alpha_1  v_{xx} = 0, &  ~~ (x,t)\in (0,L)\times \mathbb{R}^+ \\
  & v(0,t)= 0, ~ \alpha_1 v_x(L,t)= -\frac{ \gamma V(t)}{  h},&  t\in \mathbb{R}^+ \\
  &(v, \dot v)(x,0)=(v^0, v^1),&  x \in [0,L]
\end{empheq}
\end{subequations}
where $\rho, \alpha_1, \gamma$ denote mass density,  elastic stiffness, and piezoelectric coefficients of the beam, respectively, $V(t)$ denotes the voltage applied at the electrodes, and $v$ denotes the longitudinal displacement of the beam. In these models an elliptic-type differential equation for the electrical component is obtained due to Gauss' law  \ref{Maxwell}.
Solving this equation and then substituting into the mechanical equations leads to the wave equation (\ref{or}). (See (\ref{homo-vol}) with $\mu \ddot p\equiv 0$. )  The system (\ref{or}) is a well-posed boundary control problem on an appropriate Sobolev space. As a side note, both Kirchhoff and Mindlin-Timoshenko small displacement assumptions yield the same stretching equations (\ref{or}).
%\begin{subequations}
%  \label{or2}
%\begin{empheq}[left={\phantomword[r]{0}{ }  \empheqlbrace}]{align}
% &  v_{tt}- D \mc L v = 0   & \quad (x,t)\in \Omega \times \mathbb{R}^+ \\
%  & \left.v \right|_{\Gamma_0}=0, \quad \left. D \mc \mc B \right|_{\Gamma_1}= -\frac{ e V(t)}{\varepsilon h}& \quad  t\in \mathbb{R}^+\\
%  &(v, \dot v)(x,0)=(v^0, v^1)& \quad  x \in \Gamma=\Gamma_0\cup \Gamma_1=\partial \Omega
%\end{empheq}
%\end{subequations}
%where $\Gamma_0\cap \Gamma_1=\emptyset,$ $v=(v_1, v_2)^{\rm T}$ denotes the longitudinal displacements in $x_1$ and $x_2$ directions,   $D$ denotes the flexural rigidity, and $\mc L=((L(v))_1, (L(v))_1)$ and $\mc B$ denote the Lam\'{e} operator and the corresponding boundary operator, respectively (see \cite{Lagnese-Lions}).
From the control theory point of view, it is well-known that a single wave equation (\ref{or}) can be exactly controlled in the energy space (therefore the uncontrolled system, i.e. $V(t)\equiv 0,$ is exactly observable). With a mechanical feedback controller in the form of  boundary damping $V(t)=v_t(L,t),$  the solutions of the closed-loop system are exponentially stable in the energy space (i.e. \cite{Komornik-P} and references therein).

%As a side note, The in-plane displacements for the analogous plate problem is modeled by the well-known Lam\'{e} system (See \cite{Lagnese-Lions}).  The Lam\'{e} system  is also well-known that the Lam\'{e} system is controllable and stabilizable through a portion of its boundary.

Exact observability and exponential stabilizability if magnetic effects are included in the mathematical models is investigated in this paper. In the fully dynamic approach,  magnetic effects are included, and hence the wave behavior of the electromagnetic fields.   We obtain a strongly coupled system of wave equations, one for stretching and one for magnetic effects.  Voltage control comes into the play through only one boundary condition at one end.  The problem of exponentially stabilizability is essentially one of simultaneous stabilizability since a single control needs to stabilize two coupled wave equations. Simultaneous control problems for wave and beam systems have been studied by a number of researchers, including \cite{Zuazua}, \cite{Komornik-P}, \cite{Lions1}, \cite{Russell}, \cite{Weiss-Tucsnak1}.  In \cite{Weiss-Tucsnak1} conditions for simultaneous exact controllability are obtained for decoupled  systems with the same input function.  In \cite{Komornik-P} the controllability of coupled strings  of different lengths connected at one end point is considered. It is shown that controllability in finite time, in  a smaller  space than the natural energy space, is determined by the ratios of the string lengths.  Simultaneous controllability for general networks are considered in \cite{Zuazua}.

It is proven here that for almost all choices of
system parameters a simple electrical feedback controller (current flowing through the electrodes) yields strong
stability. However, for almost all system parameters, the uncontrolled system is not exactly observable in the energy space, and therefore there is no feedback $V(t)$ that makes the system exponentially stabilizable in the energy space. Finally, it is shown that the system can be exponentially stabilized only  for a  set of system parameters of Lebesque measure zero.
This behavior is qualitatively very different from the electrostatic or quasi-static models.
%\kchg{ Ozkan discuss   [D], [JT], [MT1], [MT2] mentioned by reviewer here? Should be in intro. I didn't realize that there wasn't something. }

This paper is organized as follows. In Section \ref{Sec-II},  a variational approach is used to derive the model; a system of partial differential equations  that include  magnetic effects. In section \ref{Sec-III}, well-posedness of the model is shown and also and strong stabilizability for a class of parameters.  Strong stabilizability is achieved with a feedback operator that is dual to the control operator. This feedback is purely electrical. Finally in Section \ref{Sec-IV}  observability and exponential stabilizability is shown to depend on    system parameters. If the system is exponentially stabilizable,  exponential stability  is achieved with  the same electrical feedback.
%  We show that
%\begin{description}
%  \item[(i)] For almost all choices of system parameters,  the closed-loop system  (\ref{pdes-stabil}) (voltage control-current measurement) is strongly stable.
%  \item[(ii)] For almost all choices of parameters, there is no feedback $V(t)$ such that the system (\ref{homo-vol}) can be exponentially stabilized in the energy space.
 % \item [(iii)] For a small set of system parameters of Lebesgue measure zero, the closed-loop system (\ref{pdes-stabil}) is exponentially stable in the energy space.
%\end{description}

Throughout this paper, dots ($\dot{w}$) indicates differentiation with respect to time, $\frac{d}{d x_1}= \frac{d}{d x}$ and $\partial \Omega$ indicates the boundary of  the beam $\Omega.$ % the electroded region and the insulated region.

 \begin{table}[h]

 {\footnotesize{
\renewcommand{\arraystretch}{1.3}
   \begin{tabular}[c]{|c|l|c|l|}
\hline
$A$ & Magnetic potential vector & $i_b$ & Volume current density\\
  \hline
  $B$ &  Magnetic flux density vector  &   $n$ & Surface unit outward normal vector \\
    \hline
  $\beta$ &  Impermittivity coefficients &  $\sigma_s$ & Surface charge density \\
    \hline
  $c,$ $\alpha$ &  Elastic stiffness coefficients & $\sigma_b$ & Volume charge density \\
    \hline
  $\gamma$ &  Piezoelectric coefficients &  $S$ & Strain tensor  \\
    \hline
$D$ &  Electric displacement vector & $T$ & Stress tensor \\
    \hline
  $E$ & Electric field intensity vector   &   $v$ & Longitudinal displacement  \\
    \hline
  $\varepsilon$ &  Permittivity coefficients &  $h$ & Thickness of the beam\\
  \hline
  $f_1$ & Lateral force resultant in $x_1$ direction &  $H$ & Magnetic field intensity vector \\
  \hline
  $\tilde f_1$ & Lateral force in $x_1$ direction & $V$ & Voltage (constant in space)\\
    \hline
  $f_3$ & Transverse force resultant in $x_3$ direction & $w$ & Transverse displacement\\
    \hline
  $\tilde f_3$ & Transverse force in $x_3$ direction &  $\mu$ & Magnetic permeability of beam\\
    \hline
 $U_1$ & $x_1$ component of the displacement field  &  $\rho$ & Mass density per unit volume\\
    \hline
$U_3$ & $x_3$ component of the displacement field & $\phi$ & Electric potential\\
    \hline
  $i_s$ & Surface  current density & &\\
    \hline
\end{tabular}
}}

\caption{Notation}
\end{table}

\section{Piezoelectric beam model with magnetic effects}

%\subsection*{Mechanical beam model}
\label{Sec-II}
Let $x_1, x_3$ be the longitudinal and transverse directions, respectively. Let the piezoelectric beam occupy the region $\Omega=[0,L]\times [-\frac{h}{2}, \frac{h}{2}]$  where $h << L .$ A very widely-used linear constitutive relationship \cite{Tiersten} for piezoelectric beams  is
\begin{eqnarray}
\label{cons-eqq10}
\left( \begin{array}{l}
 T \\
 D \\
 \end{array} \right)=
\left[ {\begin{array}{*{20}c}
   c & -\gamma^{\text{T}}  \\
   \gamma & \varepsilon  \\
\end{array}} \right]\left( \begin{array}{l}
 S \\
 E \\
 \end{array} \right)
\end{eqnarray}
where $T=(T_{11}, T_{22}, T_{33}, T_{23}, T_{13}, T_{12})^{\text T}$ is the stress vector, \\$S=(S_{11}, S_{22}, S_{33}, S_{23}, S_{13}, S_{12})^{\text T}$  is the strain vector, $D=(D_1, D_2, D_3)^{\text T}$ and  $E=(E_1, E_2, E_3)^{\text{T}}$ are the electric displacement  and the electric field vectors, respectively, and moreover, the matrices $[c], [\gamma], [\varepsilon]$ are the matrices with elastic, electro-mechanic and dielectric constant entries (for more details the reader can refer to \cite{Tiersten}). A list of all notation used for the piezoelectric beam model is in Table 1.   Under the assumption of transverse isotropy and polarization in $x_3-$direction, these matrices reduce to
\begin{eqnarray}
\nonumber &c= \left[ {\begin{array}{*{20}c}
   {c_{11} } & {c_{12} } & {c_{13} } & 0 & 0 & 0  \\
   {c_{21} } & {c_{22} } & {c_{23} } & 0 & 0 & 0  \\
   {c_{31} } & {c_{32} } & {c_{33} } & 0 & 0 & 0  \\
   0 & 0 & 0 &    c_{44} & 0 & 0 \\
      0 & 0 & 0 &    0 & c_{55} & 0 \\
         0 & 0 & 0 &    0 & 0 & c_{66}
\end{array}} \right],~~~~  \gamma= \left[ {\begin{array}{*{20}c}
   0 & 0 & 0 &    0 & \gamma_{15} & 0 \\
      0 & 0 & 0 &    -\gamma_{15} & 0 & 0 \\
         \gamma_{31} & \gamma_{31} & \gamma_{33} &    0 & 0 & 0
\end{array}} \right]&\\
\nonumber  &\varepsilon= \left[ {\begin{array}{*{20}c}
   \varepsilon_{11} & 0 & 0 \\
      0 & \varepsilon_{22} & 0\\
         0 &  0 & \varepsilon_{33} \\
\end{array}} \right].&
\end{eqnarray}

Since $h << L$,  assume that all forces acting in the $x_2$ direction are zero. Moreover, $ T_{33}$ is also assumed to be zero. Therefore
$$T=(T_{11}, T_{13})^{\text T}, S=(S_{11}, S_{13})^{\text T}, D=(D_1, D_3)^{\text T}, E=(E_1, E_3)^{\text{T}}$$ and (\ref{cons-eqq10}) reduces to
\begin{eqnarray}
\nonumber
\left( \begin{array}{l}
T_{11} \\
 T_{13} \\
 D_1 \\
 D_3
 \end{array} \right)
&= \left( {\begin{array}{*{20}c}
   {c_{11} }  &  0 & 0 & -\gamma_{31}  \\
   0 & {c_{55} } & -\gamma_{15}  & 0 \\
   0 & \gamma_{15} & \varepsilon_{11}&   0 &  \\
      \gamma_{31} & 0 & 0  & \varepsilon_{33}  \\
\end{array}} \right)\left( \begin{array}{l}
  S_{11} \\
  S_{13} \\
 E_1 \\
 E_3
 \end{array} \right).
\end{eqnarray}
Finally, for an Euler-Bernoulli beam, the shear stress $ S_{13} =0$. (See \ref{shear}. ) The linear constitutive equations for an Euler-Bernoulli piezoelectric beam are thus
\begin{subequations}
\label{cons-eq10}
\begin{empheq}[left={\phantomword[r]{0}{ }  \empheqlbrace}]{align}
\label{cons-eq11}  &  T_{11}=c_{11} S_{11}-\gamma_{31}E_3 & \\
\label{cons-eq11a} &T_{13}= -\gamma_{15}E_1&\\
\label{cons-eq12} & D_1=\varepsilon_{11}E_1 &\\
 \label{cons-eq13}& D_3=\gamma_{31} S_{11}+\varepsilon_{33}E_3.
\end{empheq}
\end{subequations}

\subsection*{Lagrangian}
 Let $\mb K, \mb P, \mb E$ and  $\mb B$  denote kinetic, potential, electrical, magnetic energies of the beam, respectively, and $\mb W$ is the work done by the external forces. Moreover, $\mb P-\mb E + \mb B$ is often called electrical enthalpy.

To model charge or current-controlled piezoelectric beams, that is, charge density or current density are prescribed at the electrodes, the pair $(S,E)$ are taken to be the independent variables. The  Lagrangian \cite{Lee,O-M}
\begin{eqnarray}\label{Lag}  \mb{L}= \int_0^T \left[\mb{K}-(\mb{P}-\mb{E} + \mb B)+ \mb{W}\right]dt\end{eqnarray}
with  constitutive equations (\ref{cons-eq10}) is appropriate.
For the Lagrangian $\mb L,$ the work done by the external forces is
 $$\mb W=\int_\Omega \left(\tilde f_1 U_1 + \tilde f_3 U_3\right) ~dX + \int_{\partial\Omega} \bar \sigma_s\phi~ d\Gamma$$ where $\tilde f_1, \tilde f_3$ are external  lateral and transverse  forces respectively,
 $\bar \sigma_s$ is the surface charge prescribed at the electrodes, and $(U_1, U_3)$ is the displacement field (see (\ref{kirc}));  the external forces are as defined in \cite{Lagnese-Lions}.
  Therefore $$\delta {\mb W}=\int_\Omega \left(\tilde f_1 \delta U_1 + \tilde f_3 \delta U_3 \right)~dX+ \int_{\partial\Omega}  \bar \sigma_s\delta \phi ~d\Gamma.$$

  For  voltage-driven electrodes, voltage is prescribed at the boundaries, and  a different Lagrangian is needed so that the applied voltage appears in the work term.
  Applying a Legendre transformation to $\mb L$ yields
    \begin{eqnarray}\label{tildeL}  \mb{\tilde L}= \int_0^T \left[\mb{K}-(\mb{P}+\mb{E})+\mb B +\mb{W}\right]~dt\end{eqnarray}
    where $\mb P+\mb E$ is the total stored energy of the beam. The new Lagrangian $\tilde{L}$ is a function of independent variables $(S, D).$ The constitutive relationship  (\ref{cons-eq10}) transforms to the following relationship for    $(T, E) $
    \begin{subequations}
\label{cons-eqcon10}
\begin{empheq}[left={\phantomword[r]{0}{ }  \empheqlbrace}]{align}
\label{cons-eqcon11}  &  T_{11}=\alpha S_{11}-\gamma\beta D_3 & \\
\label{cons-eqcon11a}  &  T_{13}=-\gamma_{1}\beta_{1} D_1 & \\
\label{cons-eqcon12} & E_1=\beta_{1}D_1 &\\
 \label{cons-eqcon13}& E_3=-\gamma\beta S_{11}+\beta D_3
\end{empheq}
\end{subequations}
where
\begin{eqnarray}
\label{coef}\gamma=\gamma_{31},~~\gamma_1=\gamma_{15}~~ \alpha=\alpha_1 + \gamma^2 \beta, ~~\alpha_1=c_{11}, ~~\beta=\frac{1}{\ep3}, ~~\beta_{1}=\frac{1}{\varepsilon_{11}}. \end{eqnarray}

%The two Lagrangians $\mb L$ and $\mb {\tilde L}$ are different since we have changed  the set of  independent variables from  $\{S, E\}$ to  $\{S, D\}.$ As a result, the work $\mb W$ done by the external forces changes (see (\ref{ext-vol})).
Calling $\delta (\cdot)$ the variation of the corresponding quantity, $\tilde {\mb L}$ in (\ref{tildeL}) is  obtained by applying the Legendre transformation to $\mb L:$
 \begin{eqnarray}\delta  {\mb L}= \int_0^T \left(\delta \mb K-\delta\left(\mb H+\int_\Omega D_iE_i ~dX\right) +\delta\left(\mb W + \int_{\partial\Omega} \sigma_s \phi~ d\Gamma\right) \right)dt=0,
 \end{eqnarray}
 where  $\phi$ is the electric potential, $\mb H$ is the enthalpy \cite{Lee} and $$\delta \mb{H}=\int_\Omega \left(T_{ij} \delta S_{ij} -D_k\delta E_k + M\cdot \delta B \right)~ dX$$
 where  $M=\frac{1}{\mu} B$ is the magnetic flux vector and $\mu$ is the permeability of the beam. The new Lagrangian $\tilde {\mb L}$ essentially remains the same since
 \begin{eqnarray}\int_\Omega \delta (D_i E_i) ~dX  &=&-\int_{\partial\Omega} \delta (\phi  D_i n_i)~ d\Gamma + \int_\Omega \delta (\phi \nabla \cdot D )  - \int_\Omega \delta ( M \cdot B ) ~dX %-\int_\Omega \delta ( D \cdot \dot A ) ~dX \\
 %\nonumber &=&\int_{\partial\Omega} \delta (\sigma_s \phi) ~d\Gamma + \int_\Omega \delta (\dot D \cdot A ) ~dX \\
 %\nonumber &=&\int_{\partial\Omega} \delta (\sigma_s \phi) ~d\Gamma + \int_\Omega \delta ((\nabla \times M) \cdot A ) ~dX \\
 % \nonumber &=&\int_{\partial\Omega} \delta (\sigma_s \phi) ~d\Gamma - \int_\Omega \delta ( M \cdot B ) ~dX
\label{ext-char}
\end{eqnarray}
 %To see difference  between the two Lagrangians in terms of the external forces, let us look at the following case.
 However, for the Lagrangian $\tilde {\mb L},$ the work done by the external forces is given by
 \begin{eqnarray}\mb W=\int_\Omega \left(\tilde f_1 U_1 + \tilde f_3 U_3\right) ~dX + \int_{\partial\Omega} \sigma_s \bar \phi ~d\Gamma
  \label{ext-vol}
  \end{eqnarray}
  where  $\bar \phi$ (namely voltage) is the electric potential prescribed at the electrodes,
   and therefore, using (\ref{charge_boun}),
   \begin{eqnarray}\nonumber \delta {\mb W}&=&\int_\Omega \left(\tilde f_1 \delta U_1 + \tilde f_3 \delta U_3\right)~dX +
    \int_{\partial\Omega} \bar \phi~ \delta \sigma_s ~d\Gamma \\
    &=&\int_\Omega \left(\tilde f_1 \delta U_1 + \tilde f_3 \delta U_3\right)~dX -    \int_{\partial\Omega} \bar \phi ~(\delta D_i)n_i ~d\Gamma \, .
    \end{eqnarray}
  Therefore, depending on the prescribed quantity at the electrodes, Lagrangian can be chosen either $\mb L$ or $ \tilde {\mb L}.$
  In this paper,  the  voltage at the electrodes is controlled.

Returning to the linear theory of Euler-Bernoulli beam small-displacement assumptions, the displacement field is
\begin{eqnarray}
\label{kirc} && U_1=v-x_3 \frac{\partial w}{\partial x_1}, ~ U_3=w
\end{eqnarray}
where $v=v(x_1)$ and $w=w(x_1)$ denote the longitudinal displacement of the center line, and transverse    displacement of the beam, respectively. Since
\begin{eqnarray}{{S_{13}=\frac{1}{2}\left(\frac{\partial U_1}{\partial x_3} + \frac{\partial U_3}{\partial x_1}\right)=0,}}\label{shear}\end{eqnarray} then the only strain component is given by
\begin{eqnarray}S_{11}=\frac{\partial U_1}{\partial x_1}=\frac{\partial v}{\partial x_1}-x_3 \frac{\partial^2 w}{\partial x_1^2}.\label{strains}
\end{eqnarray}

 \subsection*{Magnetic effects}
  The magnetic energy  is added to the Lagrangian $\tilde{\mb L}$ through Maxwell's equations. Let $B$ denote magnetic field vector, and $\sigma_b, i_b, \sigma_s, i_s, V, \mu, n$ denote body charge density, body current density, surface charge density,
 surface current density, voltage, magnetic permeability, and unit normal vector respectively.  Maxwell's equations are
  \begin{subequations}
  \label{Maxwell}
\begin{empheq}[left={\phantomword[r]{0}{ }  \empheqlbrace}]{align}
\label{Gauss-law} \nabla\cdot D =~\sigma_b & \quad{\rm{in}} \quad\Omega \times \mathbb{R}^+~\quad&\text{(Electric Gauss's ~law)}\\
\label{Gauss-magne} \nabla\cdot \textrm{B}=~0 &\quad{\rm{in}} \quad\Omega \times \mathbb{R}^+~\quad  & \text{(Gauss's law of magnetism)}\\
\label{Faraday} \nabla\times E=~-\dot {\textrm{B}} & \quad{\rm{in}} \quad\Omega \times \mathbb{R}^+~\quad&\text{(Faraday's law)}\\
\label{Ampere}  \frac{1}{\mu}(\nabla\times \textrm{B})= ~i_b +  \dot D & \quad{\rm{in}} \quad\Omega \times \mathbb{R}^+~\quad& \text{(Amp\'{e}re-Maxwell law)}
\end{empheq}
\end{subequations}
with one of the essential electric boundary conditions prescribed on the electrodes
\begin{subequations}
\begin{empheq}[left={\phantomword[r]{0}{ }  \empheqlbrace}]{align}
\label{charge_boun}  - D\cdot n =~\sigma_s  & \quad{\rm{on}} \quad {\partial\Omega} \times \mathbb{R}^+~\quad&\text{(Charge )}\\
\label{current_boun} \frac{1}{\mu}(\textrm{B} \times n) =~i_s & \quad{\rm{on}} \quad {\partial\Omega} \times \mathbb{R}^+~\quad & \text{(Current)}\\
\label{voltage_boun} \phi=~V  & \quad{\rm{on}} \quad {\partial\Omega} \times \mathbb{R}^+~\quad& \text{(Voltage)}
\end{empheq}
\end{subequations}
and with a chosen mechanical boundary condition at the edges of the beam  (the beam is clamped, hinged, free, etc.).
Since the electrodes are  voltage-driven,  (\ref{voltage_boun}) is appropriate.

  In modeling piezoelectric beams, there are mainly three approaches to include electric and magnetic effects \cite{Tiersten}:
\begin{description}
  \item{1) Electrostatic electric field:} An electrostatic electric field  is the most widely-used approach.
  It completely ignores magnetic effects: $\textrm{B}=\dot D=i_b=\sigma_b=0.$  Maxwell's equations  (\ref{Maxwell}) reduce to
   $\nabla\cdot D =0$ and $\nabla \times E=0.$ Therefore, by Poincar\'{e}'s theorem \cite{Dautray-Lions} there exist a scalar electric potential
   such that $E=-\nabla \phi$  and $\phi$ is determined up to a constant.

  \item{2) Quasi-static electric field:} approach ignores some of the magnetic effects (polarizable but non-magnetizable materials) \cite{Tiersten}: it is allowed that $\dot D$ and $\textrm{B}$ are non-zero, however $\sigma_b=i_b=0.$  Therefore, (\ref{Maxwell}) reduces to
       \begin{eqnarray}\nonumber\nabla \cdot D=0, ~~~\nabla\cdot \textrm{B}=0, ~~~\dot {\textrm{B}}= -\nabla\times E, ~~~ \dot D= \frac{1}{\mu}(\nabla\times \textrm{B}).\end{eqnarray} The equation $\nabla\cdot \textrm{B}=0$
       implies that there exists a magnetic potential vector $A$ such that $\textrm{B}=\nabla\times A,$ by Poinc\'{a}re's theorem. It follows from substituting $\textrm{B}$ to $\dot {\textrm{B}}= -\nabla\times E$ that there exists a scalar electric potential $\phi$ such that $E=-\nabla\phi-\dot A.$
       %The magnetic potential $A$ is not unique (see \cite{O-M}).
       One simplification in this approach is to  set $A=0$ and $\dot A=0$ since $A, \dot A\ll \phi.$  Note that $\dot D$ non-zero.
  \item {3) Fully dynamic electric field:} Unlike the quasi-static assumption,  $A$ and $\dot A$ are left in the model.  Depending on the type of material, body charge density $\sigma_b$ and body current density $i_b$ can also be  non-zero. Note that even though the  displacement current $\dot D$ is assumed to be non-zero in both quasi-static and fully dynamic approaches, the term $\ddot D$  is  zero in quasi-static approach since $\dot A=0.$
\end{description}

%The electrostatic assumption has been widely used in the literature to reduce the complexity of Maxwell's equations due the minor magnetic effects on the dynamics of piezoelectric structures. %However, these effects become significant when it comes to the controller design.

In this paper,  the third, dynamic, approach is used for the modeling of a piezoelectric beam.
Assume that there is neither external body  charges nor body currents, i.e., $\sigma_b=i_b\equiv 0$.
The magnetic field $\textrm{B}$ is perpendicular to the $x_1-x_3$ plane due to (\ref{current_boun}), and therefore
$\textrm{B}$ has only the $y-$component $\textrm{B}_2,$ and it is only a function of $x_1=x.$  This is simply because the surface current
$i_s$ at the electrodes have only $x-$component (tangential)  and $\textrm{B}$ is perpendicular to both the
outward normal vector ($n=(0,0,1)$ or $n=(0,0,-1)$) at the electrodes and $i_s.$  Also assume that $E_1=0,$ and thus $D_1=0$ by (\ref{cons-eqcon12}). Therefore, Maxwell's equations including the effects of $\textrm{B}$ become
$$ \nabla\cdot \textrm{B}=0, ~~~\dot {\textrm{B}}= -\nabla\times E, ~~~ \dot D= \frac{1}{\mu}(\nabla\times \textrm{B}).$$
It follows from the last equation that  $\frac{ d \textrm{B}_2}{dx}=-\mu \dot D_3,$ and so $$\textrm{B}_2=-\mu\int_0^x \dot D_3(\xi, x_3, t) ~d\xi.$$
The magnetic energy, which can be regarded as the   ``electric kinetic energy'', is
\begin{eqnarray*}
\mb B&=& \frac{1}{2\mu}\int_{\Omega} \|\textrm{B}\|^2~dX =  \frac{1}{2\mu}\int_{\Omega} (\textrm{B}_2)^2~dX =   \frac{\mu}{2}\int_{\Omega} \left[\int_0^x \dot D_3(\xi, x_3, t) d\xi\right]^2 ~dX.
\end{eqnarray*}
The next assumption is that $D_3$ does not vary in the thickness direction $$D_3(x,x_3,t)=D_3(x,t).$$ This assumption lines up with choice of electrical potential $\varphi(x,z,t)$ defined above to be linear in the thickness direction \cite{Rogacheva}, i.e. $\varphi(x,z,t)=\varphi^0(x,t) + z\varphi^1(x,t).$ Therefore the electric field component in the thickness direction satisfies $$E_3=\frac{\partial \varphi}{\partial z}= \beta D_3(x,t).$$
%*OZKAN: Explain previous sentence and give a physical interpretaion of $p$: electrical potential in one direction?*
Now define \begin{eqnarray}\label{defp}p=\int_0^x  D_3(\xi,  t) ~d\xi\end{eqnarray} to be the total electric charge at point $x.$ Therefore $p_x=D_3.$

\subsection*{Hamilton's Principle}
 Using  (\ref{cons-eqcon10}) (with $D_1=0$), (\ref{cons-eq10}), (\ref{strains}), and the definition (\ref{defp}) of $p$,  the stored energy (potential+ electric) $\mb P+\mb E$, magnetic  energy $\mb B$  and kinetic energy $\mb K$ of the beam are \begin{eqnarray}
\nonumber \mb P+\mb E &=& \frac{1}{2}\int_\Omega \left(T_{11}S_{11} + D_3 E_3\right) ~dX\\
%  \int_{-h/2}^{h/2} \left(\alpha_1   (v_x-x_3 w_{xx})^2+ E_3 (v_x-x_3 w_{xx})\right)dx_3\\
\label{stored-energy} &=& \frac{h}{2} \int_0^L \left[\alpha  \left(v_x^2 +\frac{h^2}{12}w_{xx}^2\right) -2\gamma\beta  v_x p_x + \beta p_x^2\right]~dx,\\
\label{m-energy} \mb{B}&=& \frac{1}{2\mu}\int_\Omega \|\textrm{B}\|^2 ~dX= \frac{\mu h}{2}\int_{0}^L \dot p^2~dx,\\
\label{k-energy} \mb{K}&=& \frac{\rho}{2} \int_\Omega \left(\dot U_1^2+ \dot U_3^2\right)~dX= \frac{\rho h}{2} \int_0^L \left[\dot v^2+\frac{ h^2}{12} \dot w_x^2 + \dot w^2\right]~dx .
\end{eqnarray}
Defining
$$f_1(x,t)=\int_{-h/2}^{h/2} \tilde f_1(x,z,t)~dz, \quad f_3(x,t )=\int_{-h/2}^{h/2} \tilde f_3(x,z,t)~dz$$
to be the external force resultants defined as in \cite{Lagnese-Lions},  and $V(t)$ the voltage applied at the electrodes,
the work done by the external forces is
\begin{eqnarray*}
\label{work-done} \mb{W}&=& \int_\Omega \left(\tilde f_1 U_1+\tilde f_3 U_3 \right)~dX- \int_{\partial\Omega} D_3 ~\bar\phi ~d\Gamma \\
&=& \int_0^L \left(f_1 v+  f_3 w- p_x V(t) \right)~dx \\
&=& \int_0^L -p_x V(t) ~ dx
\end{eqnarray*}
since there is no applied external force $\tilde{f}_1$ or lateral force $\tilde{f}_2$.

Application of Hamilton's principle, setting the variation of admissible displacements $\{v,w,p\}$ of $\tilde {\mb L}$ to zero, yields two sets of equations one for stretching and one for bending with  associated boundary conditions
\begin{subequations}
  \label{homo-vol}
\begin{empheq}[left={\phantomword[r]{0}{\text{Stretching: ~} }  \empheqlbrace}]{align}
 \label{eq-main1} &\rho  \ddot v-\alpha   v_{xx}+\gamma  \beta    p_{xx} = 0   & \\
 \label{eq-main2}  & \mu  \ddot p   -\beta   p_{xx} + \gamma \beta v_{xx}= 0, &
\end{empheq}
\begin{empheq}[left={\phantomword[l]{}{ } \empheqlbrace}]{align}
 & v(0)= p(0)=\alpha  v_{x}(L)-\gamma \beta p_x(L)=0, ~~\beta  p_x(L) -\gamma \beta v_x(L)= -\frac{V(t)}{h}& \\
  &(v, p, \dot v, \dot p)(x,0)=(v^0,  p^0,  v^1,  p^1).&
  \label{ivp}
\end{empheq}
\end{subequations}

  \begin{subequations}
  \label{i-b-c}
\begin{empheq}[left={\phantomword[r]{0}{ \text{Bending:} ~ } \empheqlbrace}]{align}
 \label{eq-main3} & \rho h \ddot w + \frac{\rho h^3}{12} \ddot{w}_{xx}  + \frac{\alpha h^3}{12}w_{xxxx}=0, &
\end{empheq}
\label{eq-main-EB-BC}
\begin{empheq}[left={\phantomword[r]{0}{ } \empheqlbrace}]{align}
\label{eq-main3-BC} & w(0)=w_x(0)=w_{xx}(L)=w_{xxx}(L)=0&\\
\label{eq-main4-BC}   &(w, \dot w)(x,0)=(w^0, w^1).&
\end{empheq}
\end{subequations}
Equation (\ref{i-b-c}) is the Rayleigh beam equation for bending. Neglecting the moment of inertia term $\frac{\rho h^3}{12}\ddot w_{xx}$ in  (\ref{i-b-c}), leads to the familiar Euler-Bernoulli beam equation. Use of  Mindlin-Timoshenko small displacement assumptions instead of Euler-Bernoulli leads to the same stretching equation (\ref{homo-vol})  \cite{accpaper}. However, the  equations for the bending and rotation of the beam are different:
  \begin{subequations}
  \label{i-b-c-MT}
\begin{empheq}[left={\phantomword[c]{}{ \quad\quad\quad\quad\quad\quad\quad\quad\quad } \empheqlbrace}]{align}
 \label{eq-main4} & \rho h \ddot w -\varsigma h (\psi + w_x)_x=0, &\\
  \label{eq-main5} & \frac{\rho h^3}{12}\ddot \psi  -\frac{\alpha h^3}{12} \psi_{xx} + \varsigma h (\psi+w_x)=0, &
\end{empheq}
\label{eq-main-MT-BC}
\begin{empheq}[left={\phantomword[l]{}{ } \empheqlbrace}]{align}
\label{eq-main4-BC-MT} & \psi (0)=\psi_x (L) = w(0)= (\psi+w_x)(L)=0&\\
\label{eq-main5-BC-MT}  &(w, \psi, \dot w, \dot \psi)(x,0)=(w^0, \psi^0, w^1, \psi^1)&
\end{empheq}
\end{subequations}
where $\psi$ and $\varsigma$ denote the angle of rotation of the beam and  shear stiffness coefficient, respectively.

%*OZKAN put correct constant in front of inertia term and add a few sentences about including other mechanical effects leads to Timoshenko*

Note that the bending equation (\ref{i-b-c}) in the Euler-Bernoulli beam case, and the bending and rotation equations (\ref{i-b-c-MT}) in  the Mindlin-Timoshenko case are completely decoupled from the stretching equations (\ref{homo-vol}).
The applied voltage $V(t)$ affects only the stretching motion. Therefore throughout the rest of the paper  only the stretching equations (\ref{homo-vol}) are considered.

Note that in the case of static magnetic effects, then $\mu \ddot p=0$ in (\ref{eq-main2}) and (\ref{eq-main2}) can be solved for $p_{xx}.$ Elimination of $p_{xx}$ in (\ref{eq-main1}) yields the system (\ref{or}).
This is the stretching equation obtained for a single piezoelectric beam in all of the classical models, i.e. \cite{Banks-Smith}, \cite{Smith}, \cite{Tiersten}. This model is known to be exactly observable and  stabilizable, i.e. see \cite{Komornik-P}. Similarly, the case of no electro-mechanical coupling, $\gamma = 0$, the voltage  $V$ only affects $p .$ We will assume throughout this paper that $\gamma >0$ and $\mu >0$ so that the stretching equations (\ref{homo-vol}) are coupled.

% We choose  $V(t)=~\frac{h}{2}\dot p(L)$ in (\ref{homo-vol}) where $\dot p(L)$ is the current flowing through the electrodes. The closed-loop system is
%\begin{subequations}
%  \label{pdes-stabil}
%\begin{empheq}[left={\phantomword[l]{}{ }  \empheqlbrace}]{align}
% \label{pdes-stabil-a} &\rho  \ddot v-\alpha   v_{xx}+\gamma  \beta    p_{xx} = 0   & \\
%\label{pdes-stabil-b} & \mu \ddot p   -\beta   p_{xx} + \gamma  \beta  v_{xx}= 0, &
%\end{empheq}
%\begin{empheq}[left={\phantomword[l]{}{ }  \empheqlbrace}]{align}
% & v(0)= p(0)= \alpha  v_{x}(L)-\gamma \beta   p_x(L)=0, ~~\beta  p_x(L) -\gamma \beta  v_x(L)= -\frac{\dot p(L)}{2}&\\
% &(v, p, \dot v, \dot p)(x,0)=(v^0, p^0, v^1, p^1).&
%\end{empheq}
%\end{subequations}

\section{Well-posedness}
\label{Sec-III}

Define
$$
  H^1_L(0,L)=\{v\in H^1(0,L): v(0)=0\}, \hspace{2em}\mX=(\Ltwo(0,L))^2  ,
$$
and  the complex linear space
$$ \mathrm{H} = \left(H^1_L(0,L)\right)^2 \times \mX .
$$
Since we are neglecting the bending terms, the  energy associated with (\ref{homo-vol})  is, recalling from (\ref{coef}) that $\alpha=\alpha_1 + \gamma^2 \beta$,
%\begin{eqnarray}
% \nonumber \mathrm{E}(t)&=&\|(v, p, \dot v,  \dot p) \|_{\mathrm E}^2\\
%\label{Energy-nat} &=&\frac{1}{2}\int_0^L \left\{\rho  |\dot v|^2 + \mu|\dot p|^2+ \alpha |v_x|^2-\gamma  v_x {\bar p}_x - \gamma {\bar v}_x p_x + %\beta |p_x|^2  \right\} dx, ~~t\in \mathbb{R}
% &=&\frac{1}{2}\int_0^L \left\{\rho  |\dot v|^2 + \mu|\dot p|^2+ \left<\left( {\begin{array}{*{20}c}
%   {\alpha } & {-\gamma }  \\
%   {-\gamma } & {\beta }  \\
%\end{array}} \right)
%\left( \begin{array}{l}
% v_x \\
% p_x \\
% \end{array} \right),
%\left( \begin{array}{l}
%v_x \\
% p_x \\
% \end{array} \right)\right>_{\mathbb{C}^2}
%\right\} dx\quad\quad\quad.
%\end{eqnarray}
 % By using (\ref{oz}), (\ref{Energy-nat}) can be written as
\begin{eqnarray}
\label{Energy-nat} \mathrm{E} &=&\frac{1}{2}\int_0^L \left\{\rho  |\dot v|^2 +   \mu |\dot p|^2 + \alpha_1   |v_x|^2 + \beta\left| \gamma v_x- p_x\right|^2  \right\}~ dx .
\end{eqnarray}
This motivates definition of the inner product on $\mathrm{H}$
{ \small{\begin{eqnarray}
\nonumber && \left<\left( \begin{array}{l}
 u_1 \\
 u_2 \\
 u_3\\
 u_4
 \end{array} \right), \left( \begin{array}{l}
 v_1 \\
 v_2 \\
 v_3\\
 v_4
 \end{array} \right)\right>_{\mathrm{H}}= \left<\left( \begin{array}{l}
 u_3\\
 u_4
 \end{array} \right), \left( \begin{array}{l}
 v_3\\
 v_4
 \end{array} \right)\right>_{(\Ltwo(0,L)^2} +  \left<\left( \begin{array}{l}
 u_1 \\
 u_2
 \end{array} \right), \left( \begin{array}{l}
 v_1 \\
 v_2
 \end{array} \right)\right>_{\left(H^1_L(0,L)\right)^2}\\
\nonumber && =\int_0^L \left\{\rho  u_3 \bar v_3 + \mu  u_4  \bar v_4\right\}~dx + \int_0^L \left\{\alpha_1   (u_1)_{x} (\bar v_1)_x\right.  \\
\nonumber  &&\quad\quad\quad\quad\quad \quad\quad\quad \left. +\beta\left( \gamma  (u_1)_x- (u_2)_x\right) \left( \gamma (\bar v_1)_x- (\bar v_2)_x\right) \right\}dx\\
\label{inner} &&= \int_0^L \left\{\rho  u_3 \bar v_3 + \mu  u_4  \bar v_4  + \left< \left( {\begin{array}{*{20}c}
   \alpha_1 + \gamma^2\beta  & -\gamma\beta\\
     -\gamma \beta & \beta  \\
\end{array}} \right) \left( \begin{array}{l}
 u_{1x} \\
 u_{2x}
 \end{array} \right), \left( \begin{array}{l}
  v_{1x} \\
  v_{2x}  \end{array} \right)\right>_{\mathbb{C}^2}\right\}~dx
 \end{eqnarray}}}
  where $\left<\cdot,\cdot\right>_{\mathbb{C}^2}$ is the inner product on $\mathbb{C}^2.$

  Rewriting the last term,
{\small{\begin{eqnarray*}
\nonumber && \left<\left( \begin{array}{l}
 v_1 \\
 v_2 \\
 v_3\\
 v_4
 \end{array} \right), \left( \begin{array}{l}
 v_1 \\
 v_2 \\
 v_3\\
 v_4
 \end{array} \right)\right>_{\mathrm{H}}=
 \int_0^L \left\{\rho  |v_3|^2 +   \mu |v_4|^2 + \alpha_1   |v_{1x}|^2 + \beta\left| \gamma v_{1x}-  v_{2x} \right|^2  \right\}~ dx,
 \end{eqnarray*}
 }}
 and so $\langle \, , \, \rangle $ does indeed define an inner product, with induced norm
 $$ \left\|\left( \begin{array}{l}
 v \\
 p \\
 \dot{v}\\
 \dot{p}
 \end{array} \right) \right\|^2  = \frac{2}{h} E.$$

Define the operator
\begin{equation}
A: {\text{Dom}}(A)\subset \mX \to \mX \quad A= \left( {\begin{array}{*{20}c}
   -\frac{\alpha }{\rho}D_x^2  & \frac{\gamma\beta }{\rho}D_x^2  \\
       \frac{\gamma\beta}{\mu} D_x^2  &  -\frac{\beta}{\mu} D_x^2  \\
\end{array}} \right) ,
\label{hom-A}
\end{equation}
where
\begin{equation}   \label{dom-hom-A} {\rm {Dom}}(A) = \{ (w_1, w_2)^{\rm T} \in (H^2(0,L)~  \cap ~ H^1_L(0,L))^2~;~ w_{1x}(L)= w_{2x}(L)= 0 \}.
\end{equation}
The operator $A$ can be easily shown to be a positive and self-adjoint operator.

For  $\theta\ge 0$  define  $\mX_{\theta}={\rm Dom} (A^{\theta})$ with the norm $\|\cdot \|_{\theta}=\|A^{\theta} \cdot\|_\mX$ The space $\mX_{-\theta}$ is  the dual of $\mX_{\theta}$ pivoted with respect to $\mX.$ For example, the inner product on $\mX_{-1/2}$ is
$$\left <z_1, z_2\right>_{\mX_{-1/2}}:=\left< A^{-1/2} z_1, A^{-1/2}z_2\right>_{\mX}.$$
Using the definition of inner product $\left<\cdot, \cdot \right>_{(H^1_L(0,L))^2}$ in (\ref{inner}) yields
\begin{eqnarray}\nonumber \left<z_1, z_2\right>_{\mX_{1/2}}&=& \left<A^{1/2} z_1, A^{1/2} z_2\right>_{\mX}=\left<A z_1, z_2\right>_{\mX}=\left< z_1, z_2\right>_{(H^1_L(0,L))^2},
\end{eqnarray}
 and therefore
  \begin{eqnarray}\label{duals}\mX_0=\mX, \quad \mX_{1/2}= (H^1_L(0,L))^2, \quad \mX_{-1/2}=((H^1_L(0,L))^*)^2   \end{eqnarray}  where $(H^1_L(0,L))^*$ is the dual space of $H^1_L(0,L)$ pivoted with respect to $\Ltwo(0,L).$ Moreover, $\mX_1= {\rm Dom}(A).$

Let $ \psi=(\psi_1, \psi_2, \psi_3, \psi_4)^{\rm T}. $
Note that $\mathrm H=\mX_{1/2}\times \mX$ and  define $\mc A: {\rm {Dom}}(\mc A) \subset \mathrm{H}\to \mathrm{H}$
by
$$ \mc A = \left( {\begin{array}{*{20}c}
   0 & I_{2\times 2}  \\
   -A & 0  \\
\end{array}} \right) ,
$$
\begin{equation}
 \begin{array}{lll}&& {\rm {Dom}}(\mc A)  = \mathrm {X}_1 \times \mathrm{X}_{1/2}\\
&& \quad\quad= \{\psi \in \mathrm H \cap ((H^2(0,L))^2\times (H^1_L(0,L))^2)  ; ~ \psi_{1x}(L)= \psi_{2x}(L)= 0 \}  \end{array}
   \label{dom-hom}
\end{equation}
which is densely defined in $\mathrm{H}.$
 Also define the control operator $B$
 \begin{eqnarray}
\nonumber  & B_0 \in \mathcal{L}(\mathbb{C}, \mX_{-1/2}), ~ \text{with} ~ B_0 =   \left( \begin{array}{c}
0 \\
-\frac{1}{h}  \delta (x-L) \end{array} \right),& \\
\label{defb_0} &\quad B \in \mathcal{L}(\mathbb{C} , \mathrm H_{-1}) , ~ \text{with} ~ B=   \left( \begin{array}{c} 0 \\ B_0 \end{array} \right)&
 \end{eqnarray}
 where $\mathrm H_{-1}$ is the dual of the space ${\rm Dom(\mc A )}=\mX_{1}\times \mX_{1/2}$ pivoted with respect to $\mathrm H=\mX_{1/2}\times \mX.$ By (\ref{duals}). We have  $\mathrm H_{-1}= \mX_0 \times \mX_{-1/2}.$
The dual operators of $B_0$ and $B$ are
\begin{eqnarray*}
& B_0^* \in \mc L( \mX_{1/2},  \mathbb{C}), ~~ B_0^* \psi = -\frac{1}{h}  \psi_4 (L),~~ \text{with} ~~  B^*\psi=(0_{2\times 2}\quad B_0^*)^{\rm T}\psi=-\frac{1}{h} \psi_4(L).&
\end{eqnarray*}
Writing $\varphi=(v, p, \dot v, \dot p)^{\rm T}$ and defining the output
$$y(t) = \frac{1}{h}\dot{p} (L,t),$$
the control  system (\ref{homo-vol})  with this output can be put into the  state-space form
\begin{subequations}
\label{Semigroup}
\begin{empheq}[left={\phantomword[l]{}{ }  \empheqlbrace}]{align}
&\dot \varphi =  \underbrace{\left( {\begin{array}{*{20}c}
   0 & I_{2\times 2}  \\
   -A & 0  \\
\end{array}} \right)}_{\mc A}  \varphi +\underbrace{ \left( \begin{array}{c} 0 \\ B_0 \end{array} \right)}_B V(t) , &\\
&\varphi(x,0) =  \varphi ^0 , \\
&y(t) = -B^* \varphi (t).
\end{empheq}
\end{subequations}

 \begin{lemma} \label{skew-adjoint}The operator $\mc{A}$  satisfies $\mc{A}^*=-\mc{A}$ on  $\mathrm{H},$ and
 \begin{eqnarray}\label{dang}{\rm Re}\left<\mc{A} \psi, \psi\right>_{\mathrm{H}}={\rm Re}\left<\mc{A}^*\psi, \psi\right>_{\mathrm{H}}= 0.\end{eqnarray}
 Also, $\mc{A} $ has a compact resolvent.
\end{lemma}

\textbf{Proof:}
%\kchg{ changing U to u and V to v and did some rewording}
Choose any $u=[u_1, u_2, u_3, u_4]^{\rm T},~v=[v_1, v_2, v_3, v_4]^{\rm T}\in \text{Dom}(\mc A).$ A simple calculation using integration by parts and the boundary conditions (\ref{dom-hom}) shows \begin{eqnarray}
 \nonumber \left<\mc{A}u , v\right>_{\mathrm{H}}  &=& \int_0^L \left\{ (-\alpha (\bar v_1)_{xx}+\gamma  \beta  (\bar v_2)_{xx})  u_3 + (\beta (\bar v_2)_{xx}-\gamma  \beta  (\bar v_1)_{xx}) u_4 \right. \\
\nonumber && ~~\left. - \alpha  (\bar v_3)_x (u_1)_x   + \gamma  \beta  (\bar v_4)_x(u_1)_{x}  +\gamma  \beta  (u_2)_{x} (\bar v_3)_{x}  - \beta  (u_2)_{x} (\bar v_4)_{x} \right\}dx\\
\label{eq110} &=& \left<u, -\mc{A} v\right>_{\mathrm{H}} . \\ %&= & \left<u. \mc{A}^* v\right>_{\mathrm{H}} .
\end{eqnarray}
This shows that $\mc A$ is skew-symmetric.
To prove that $\mc A$ is skew-adjoint on $\mathrm H,$ i.e. $\mc A^*=-\mc A$ on $\mathrm H,$ with the same domains  it is required to show that  for any  $v\in \mathrm{H}$  there is $u \in  \text{Dom}(\mc A )$ so that   $\mc A u=v. $
This is equivalent to solving the  system of equations  for $ u \in  \text{Dom}(\mc A ).$
Using (\ref{coef}) to simplify the equations leads to
   \begin{eqnarray}
 \nonumber u_3&=& v_1\\
 \nonumber u_4&=& v_2\\
%\nonumber-\alpha (u_1)_{xx} +\gamma\beta (u_2)_{xx}&=& \rho v_3\\
%\label{salak} \gamma \beta (u_1)_{xx}-\beta (u_2)_{xx}&=& \mu v_4.
\nonumber -(u_1)_{xx} &=& \frac{\rho}{\alpha_1}v_3+ \frac{\mu \gamma }{\alpha_1} v_4\\
\label{salak} -(u_2)_{xx}&=& -\frac{(\alpha+\alpha_1)\rho}{\alpha_1 \gamma \beta}v_3 - \frac{\alpha\mu}{\beta\alpha_1}v_4 \, .
 \end{eqnarray}
  % \begin{eqnarray}
% \nonumber u_3&=& y_1\\
 %\nonumber u_4&=& y_2\\
%\nonumber -(u_1)_{xx} &=& \frac{\rho}{\alpha_1}y_3+ \frac{\mu \gamma }{\alpha_1} y_4\\
%\nonumber -(u_2)_{xx}&=& -\frac{(\alpha+\alpha_1)\rho}{\alpha_1 \gamma \beta}y_3 - \frac{\alpha\mu}{\beta\alpha_1}y_4.
 % \end{eqnarray}
  Since the Greens function corresponding to the operator $-\frac{d^2}{dx^2}$ with the boundary conditions $(\cdot)(0)=\frac{d(\cdot)}{dx}(L)=0$ is
$K(x,r)=\left\{
                 \begin{array}{ll}
                   r, &  x>r \\
                    x, & x<r,
                 \end{array}
               \right.
$ the solution of (\ref{salak}) is
\begin{eqnarray}&&\nonumber u_1=\frac{1}{\alpha_1}\int_0^L K(x,r) \left(\rho v_3(r)+ \mu \gamma  v_4(r)\right)~dr\\
\nonumber &&u_2=-\frac{1}{\alpha_1}\int_0^L K(x,r) \left(\frac{(\alpha+\alpha_1)\rho}{ \gamma \beta}v_3(r) +\frac{\alpha\mu}{\beta}v_4(r)\right)~dr \\
\label{salak2} && u_3=v_1, ~~u_4=v_2.
\end{eqnarray}
Using $v\in\mathrm{H},$ i.e. $v_1, v_2 \in H^1_L(0,L)$ and $v_3, v_4\in L^2(0,L),$  implies that $u_3, u_4  \in H^1_L(0,L)$ and $u_1, u_2 \in H^2(0,L)\cap H^1_L(0,L)$  with $(u_1)_x(L)=(u_2)_x(L)=0.$ Therefore, $u\in \text{Dom}(\mc A )$ is uniquely defined. Using Proposition 3.7.3 in \cite{Weiss-Tucsnak}  leads to the conclusion that $\mc A^*=-\mc A$ on $\mathrm H.$
Since then for $u \in \text{Dom}(\mc A )$,
%$$ \left<\mc A u , u \right>= -\left< u , \mc A u \right>= -\overline{\left<\mc A u , u \right>$$
with a similar expression for $A^*$, (\ref{dang}) follows.

 Moreover, $\text{Dom}(\mc A )$ is densely
defined and compact in $\mathrm{H}$ by Sobolev's embedding theorem. Therefore, for any $\lambda \in \rho (\mc A )$,  $(\lambda I -\mc A )^{-1}$ is a compact operator. $\square$

%*OZKAN - need a proof of compact resolvent for $\mc A$ so its eigenvectors are a basis*

% $$ y(t)=-B^*\psi = (0_{2\times 2}\quad B_0^*)~ \psi $$

%\begin{eqnarray}
%\nonumber  \int_0^T |\dot p(L)|^2~dt  &\le& \frac{1}{m}\int_0^T \left(\frac{\rho L}{2}|\dot v(L)|^2 +\frac{L\mu }{2} |\dot p(L)|^2 \right)~dt \le %C(T)\mathrm{E}(t).\quad\square
%\end{eqnarray}

%The boundary control system (\ref{ivp }) is written in state space form in (\ref{\Semigroup}) and has transfer function $\mb G (s) . $
The transfer function  corresponding to the  control system (\ref{Semigroup})  is (see \cite{Weiss-a} for the calculation for a similar system)
\begin{eqnarray}\label{transfer}\mb G( s ) = s B_0^*( s ^2 I + A)^{-1}B_0
\end{eqnarray}
for $s$, $\Re s >0.$
 %Note also that The boundary conditions at $x=L$ for $u_1$ and $u_3$ in the domain (\ref{dom-hom}) imply that
%$(u_1)_x(L,t)=(u_3)_x(L,t)=0$ by using (\ref{coef}).

%vspace{0.1in}

\begin{lemma} \label{appen-1}   Define the set  $\mc C_{s_1}=\{ s \in \mathbb{C}~:~ s = s _1+ i s _2, \quad s_1 >0 \}. $  We have
\begin{eqnarray}\label{sok-1}\mathop {\sup }\limits_{ s \in \mc C_{ s_1}}\|\mb G( s )\|_{\mc L(\mathbb{C})}<\infty.\end{eqnarray}
%$A$ and $B_0$ are defined by (\ref{defb_0}). Moreover,
\end{lemma}
\textbf{Proof:} See Appendix \ref{appenB}.

\begin{definition} The operator $B \in \mc L(\mathbb{C},\mathrm H_{-1})$ is an admissible control operator for $\{e^{\mc A t}\}_{t\ge 0}$ if
there exists a positive constant $c(T)$ such that for all $u \in H^1 (0,T )$,
$$\left\|  \int_0^T  e^{\mc A (T- t)} B u (t) dt \right\|_{\mathrm H} \le c(T) \| u \|_{\Ltwo (0,T)}. $$
\end{definition}
%vspace{0.1in}

\begin{definition} The operator $B^*\in \mc L({\rm Dom}(\mc A), \mathbb{C})$ is an admissible observation operator for $\{e^{\mc A^*} t\}_{t\ge 0}$ if
there exists a positive constant $c(T)$ such that for all $\varphi^0\in {\rm Dom}(\mc A)$
$$\int_0^T \|B^* e^{\mc A^* t} \varphi^0\|^2~ dt \le c(T) \|\varphi^0\|^2_{\mathrm H}.$$
\end{definition}
%vspace{0.1in}

The operator $B^*$ is an admissible observation operator for $\{e^{\mc A^*t} \}_{t\ge 0},$ if and only if $B$ is an admissible control operator for $\{e^{\mc A t}\}_{t\ge 0}$ \cite[pg. 127]{Weiss-Tucsnak}).
%vspace{0.1in}

Consider the uncontrolled system
\begin{equation}
\label{Semigroup-H}
\begin{array}{lll}
\dot{\varphi }(t) &=& \mc A \varphi (t), \\
  \varphi(x,0)& =& \varphi ^0,  \\
y(t) &=&-B^*\varphi .
\end{array}
\end{equation}

%vspace{0.1in}

%vspace{0.1in}
The following theorem on well-posedness of (\ref{Semigroup}) and (\ref{Semigroup-H}) is now immediate.
 It proves that for any $T>0,$ the map from the input $V(t) \in \Ltwo(0,T)$ to the solution $\psi\in \mathrm {H},$ and the map from the input $V(t)$ to the output $y(t)$  of (\ref{Semigroup}) are bounded.

%vspace{0.1in}
\begin{theorem}\label{w-pf}
Let $T>0,$ and $V(t)\in \Ltwo(0,T).$ For any $\varphi^0 \in \mathrm{H},$ there exists  positive constants $c_1(T),c_2(T)$
%and a unique solution  to (\ref{Semigroup})   $\varphi (T)  \in \mathrm{H} ,$
such that
      \begin{eqnarray}\label{conc}\|\varphi (T) \|^2_{\mathrm{H}} &\le& c_1 (T)\left\{\|\varphi^0\|^2_{\mathrm{H}} + \|V\|^2_{\Ltwo(0,T)}\right\},\\
      \label{conc-a} \|y\|^2_{\Ltwo(0,T)} &\le& c_2(T) \left\{\|\varphi^0 \|_{\mathrm H}^2 + \|V\|^2_{\Ltwo(0,T)}\right\}.
      \end{eqnarray}
\end{theorem}

%vspace{0.1in}
\textbf{Proof:}  The operator $B^*$ defined above is an admissible observation operator for the system (\ref{Semigroup-H})
by  Lemma \ref{appen-1} (see Proposition 3.2 and 3.3 in \cite{Ammari-Tucsnak}). Therefore $ B$ is an admissible control operator for the semigroup $\{e^{\mc A t}\}_{t\ge 0}$ corresponding to (\ref{Semigroup}).
%We know that the transfer function $\mb G( s )$ of  (\ref{Semigroup}) is bounded on every vertical strip in the open right half plane $\{ s \in \mathbb{C}~:~ s = s _1+ i s _2,~ s _1>0\}$ by Lemma \ref{appen-1}. This is equivalent to showing that  $\mb G( s )$ is bounded on the region $\{ s \in \mathbb{C}~:~ s = s _1+ i s _2,~ s _1\ge \alpha\}$ for any $\alpha\in\mathbb{R}$ by \cite[Remark 1]{Guo-Luo}.
Lemma  \ref{appen-1} and the Paley-Wiener Theorem implies that the map from the input $V$ to the output $y$ is bounded from $ \Ltwo(0,T)$ to $\Ltwo(0,T)$ \cite[Thm. 5.1]{C-W}.
% Therefore by Theorem 5.1 in \cite {C-W}, the triple $(\mc A, B, B^*)$ is well-posed.
 The conclusions (\ref{conc}) and (\ref{conc-a}) follow. $\square$
%vspace{0.1in}

Alternatively, the state could be defined as
$$(\sqrt{\rho} v_t, \sqrt{\alpha_1}v_x, \sqrt{\mu} p_t, \sqrt{\beta} (p_x-\gamma v_x)). $$
With this choice of state, the control system is well-posed on $[\Ltwo (0,L)]^4 $ and is  a port-Hamiltonian system \cite{cdcpaper}.

%%%%%%
%% DAMPED SYSTEM
%%%%%%%%

\subsection*{Damped system}

%When we choose a feedback controller $V(t)=\frac{h}{2}\dot p(L)$ we obtain  (\ref{pdes-stabil}).
Setting the  control signal in (\ref{Semigroup}) to be $ V(t) = -\frac{1}{2} B^*z  + u(t)$  where $u(t)$ is a new controlled input and  modifying the output slightly leads to the system
\begin{subequations}
\label{Sg}
\begin{empheq}[left={\phantomword[r]{0}{ }  \empheqlbrace}]{align}
&\dot z(t) = \mc A_d z(t)  + B u(t)= \left( {\begin{array}{*{20}c}
   0 & I_{2\times 2}  \\
   -A & -\frac{1}{2} B_0 B_0^*  \\
\end{array}} \right)z + \left( \begin{array}{c}
0_{2\times 2}\\
B_0 \end{array} \right) u(t), &\\
&z(x,0) =  z^0,&\\
\label{out}&y(t)=-B^*z(t) + u(t)&
\end{empheq}
\end{subequations}
where $\mc A_d : \text{Dom} (\mc A_d ) \subset \mathrm H \to \mathrm H$ and $\text{Dom} (\mc A_d )$ is defined by
 \begin{eqnarray}\nonumber &{\rm {Dom}}(\mc A_d)=\left\{z \in (H^2(0,L))^2\times (H^1_L(0,L))^2~:~ z_1(0)=z_2(0)=0,\right.&~~\\
\label{sem-dom} &\left.\alpha  z_{1x}(L)-\gamma \beta z_{2x}(L)=0, ~~\beta  z_{2x}(L) -\gamma \beta z_{1x}(L)= -\frac{z_4(L)}{2h^2} \right\}.&
\end{eqnarray}
%vspace{0.1in}
This system can also be written in second-order form as
\begin{subequations}
\label{Sg-aux}
\begin{empheq}[left={\phantomword[r]{0}{ }  \empheqlbrace}]{align}
&\left( \begin{array}{l}
 \ddot v \\
 \ddot p
 \end{array} \right)+ A \left( \begin{array}{l}
  v \\
  p
 \end{array} \right) + \frac{1}{2} B_0 B_0^* \left( \begin{array}{l}
 \dot v \\
 \dot p
 \end{array} \right)   = B_0u(t), &\\
&\left( \begin{array}{l}
  v \\
  p
 \end{array} \right)(x,0) =  \left( \begin{array}{l}
  v^0 \\
  p^0
 \end{array} \right), \quad \left( \begin{array}{l}
 \dot  v \\
 \dot  p
 \end{array} \right)(x,0) =  \left( \begin{array}{l}
  v^1 \\
  p^1
 \end{array} \right)&\\
\label{out-aux}&y(t)=-B_0^* \left( \begin{array}{l}
 \dot  v \\
 \dot  p
 \end{array} \right) + u(t). &
\end{empheq}
\end{subequations}
This system is  a member of the class studied in \cite{Weiss-a}.
%vspace{0.1in}

%*OZKAN couldn't the following just be theorem quoted from \cite{Weiss}?*
Let $\mathrm H^d_{-1}$ is the dual of the space ${\rm Dom(\mc A_d )}$ pivoted with respect to $\mathrm H=\mX_{1/2}\times \mX.$
%vspace{0.1in}
\begin{theorem} \label{eigens} Let $T>0.$ The system (\ref{Sg}) defines a well-posed and conservative linear system with the input $u(t)\in \Ltwo(0,T),$ the output $y(t)\in \Ltwo(0,T),$ the state space $\mathrm H,$ the semigroup $\{e^{\mc A_d t}\}_{t\ge 0},$ and the transfer function $\mb G_d$. Then $\mc A_d $ is the generator of a contraction semigroup on $\mc H$, $B \in \mc L(\mathbb C, \mathrm{H}^d_{-1})$ and $B^*\in \mc L( {\rm Dom}(\mc A_d ),  \mathbb C)$ are admissible control and observation operators, respectively, and $\|\mb G_d( s )\|\le 1$ for all $ s  \in C_{ s }=\{ s \in \mathbb{C}~:~ s = s _1+ i s _2, ~~~ s _1>0\}.$
\end{theorem}
%vspace{0.1in}

\textbf{Proof:}  Since $E(t)=\frac{1}{2}\|z(t)\|^2_{\mathrm H}$ by (\ref{Energy-nat}), a direct calculation by using (\ref{Sg})  reads
\begin{eqnarray}\nonumber E(T)-E(0)&=&\int_0^T \left(-\frac{1}{2}\left<B^* z, B^* z\right>_{\mathbb{C}^2} + \frac{1}{2}\left<u(t), B^*z \right>_{\mathbb{C}^2} + \frac{1}{2}\left<B^*  z, u(t) \right>_{\mathbb C^2}\right)~dx\\
\nonumber &=& \frac{1}{2}\left(\int_0^T |u(t)|^2~ dt - \int_0^T |y(t)|^2~ dt\right),
\end{eqnarray}
 and therefore
  \begin{eqnarray}\label{sok2}\|z(t)\|^2_{\mathrm H} + \int_0^T |y(t)|^2 = \|z^0\|^2_{\mathrm H} + \int_0^T |u(t)|^2.\end{eqnarray}
By Proposition 4.5 in \cite{Weiss-a}, the conclusion of the theorem follows. $\square$
 %and the feedback $B^*z =\dot p(L)=\int_0^L  \dot D_3(v,  t)~ dv$ is the  current flowing through the electrodes of the beam.

%\subsection{Strong stability of the semigroup $\{e^{\mc A_d t}\}_{t\ge 0}$}

We now show that the semigroup $\{e^{\mc A_d t}\}_{t\ge 0}$ is strongly stable for almost all choices of system parameters.
%vspace{0.1in}

\begin{theorem}\label{iso} The spectrum $\sigma({\mc A_d })$ of $\mc A_d $ has all isolated eigenvalues,  and $0\in \sigma(\mc A_d ).$
\end{theorem}
%vspace{0.1in}

\textbf{Proof:} First show that $0\in \rho(\mc A_d ).$ Let $G=(g_1, g_2, g_3, g_4)\in \mathrm{H}$ and  find $U=(u_1, u_2, u_3, u_4)$ such that $U\in {\rm{Dom}}(\mc A_d )$ and $ {\mc A_d } U=G.$
%\begin{eqnarray}
%\nonumber u_3 &=& f_1\\
%\nonumber u_4&=& f_2\\
%\nonumber \alpha (u_1)_{xx}-\gamma (u_2)_{xx}&=& \rho f_3\\
%\label{eq30} \beta (u_2)_{xx}-\gamma  (u_1)_{xx} &=& \frac{1}{\mu} f_4.
%\end{eqnarray}
Similar to (\ref{salak2}), the solution of $\mc A_d  U=G$ is
\begin{eqnarray}
\nonumber u_3 &=& g_1\\
\nonumber u_4&=& g_2\\
\nonumber u_1&=& \frac{1}{\alpha_1  } \int_0^L\left(\rho g_3(r)+\gamma \mu g_4(r)\right)K(x,r)~dr - \frac{\gamma}{2h^2\alpha_1}g_2(L)x\\
\nonumber u_2&=& -\frac{1}{\alpha_1  } \int_0^L\left(\frac{(\alpha+\alpha_1)\rho}{ \gamma \beta}g_3(r)+\frac{\mu\alpha}{\beta} g_4(r)\right)K(x,r)~dr - \frac{1}{2h^2}\left(\frac{\gamma ^2}{\alpha_1} + \frac{1}{\beta}\right)g_2(L)x
\end{eqnarray}
where $K(x,r)=
\left\{ \begin{array}{l}
 x,\quad x\le r \\
 r,\quad x\ge r. \\
 \end{array} \right.$  Since $G\in\mathrm{H},$ $g_1, g_2 \in H^1_L(0,L)$ and $g_3, g_4\in \Ltwo(0,L),$ and by the Trace theorem $g_2(L)\in \Ltwo(0,L).$ Note that $u_1$ and $u_2$ satisfy the boundary conditions in (\ref{sem-dom}). Therefore $U\in {\rm{Dom}(\mc A_d )}.$ Also, there is a unique solution $U .$ Thus $0 \in \rho (\mc A_d ) .$

 Moreover, $\text{Dom}(\mc A_d )$ is densely
defined and compact in $\mathrm{H}$ by Sobolev's embedding theorem. This together with $0\in \rho(\mc A_d ) $ implies that $(\lambda I -\mc A_d )^{-1}$
is compact at $\lambda=0,$ thus compact for all $\lambda\in \rho(\mc A_d ).$ Hence the spectrum of $\mc A_d $ contains all isolated eigenvalues. $\square$
%vspace{0.1in}

\begin{theorem} \label{stronglystable} $\{e^{\mc A_d t}\}_{t\ge 0}$ is strongly stable in $\mathrm{H}$ if and only if   $\frac{\zeta_1}{\zeta_2}\ne \frac{2n-1}{2m-1},$ for some $ n,m \in \mathbb{N}$  where
\begin{eqnarray}
\label{lam1} \zeta_1 &=& \frac{1}{\sqrt{2}}\sqrt{\frac{\gamma^2\mu }{\alpha_1}+\frac{\mu}{\beta}+\frac{\rho}{\alpha_1  }+\sqrt{\left(\frac{\gamma^2\mu}{\alpha_1 }+\frac{\mu}{\beta}
+\frac{\rho}{\alpha_1  }\right)^2-\frac{4\rho  \mu}{\beta \alpha_1 }}} \\
\label{lam2} \zeta_2 &=& \frac{1}{\sqrt{2}}\sqrt{\frac{\gamma ^2\mu }{\alpha_1}+\frac{\mu}{\beta}+\frac{\rho}{\alpha_1  }-\sqrt{\left(\frac{\gamma^2 \mu}{\alpha_1  }+\frac{\mu}{\beta}
+\frac{\rho}{\alpha_1  }\right)^2-\frac{4\rho  \mu}{\beta\alpha_1 }}}.
\end{eqnarray}
\end{theorem}
%vspace{0.1in}

\textbf{Proof:} By Theorems \ref{eigens} and  \ref{iso}, the spectrum consists of only eigenvalues, and ${\rm{Re}}\lambda\le 0.$  The eigenvalue problem
$$\mc A_d z = \lambda z $$
with $z=(v, p, \tilde{v}, \tilde{p} )$
can be written
\begin{subequations}
  \label{pdes-stabil10}
\begin{empheq}[left={\phantomword[r]{0}{ }  \empheqlbrace}]{align}
\label{st-010}  &\alpha   v_{xx}-\gamma  \beta p_{xx} = \rho\lambda^2 v   & \\
\label{st-020} &   \beta   p_{xx} -\gamma \beta v_{xx}= \mu \lambda^2 p , &\\
&\tilde{v} = \lambda v, & \\
& \tilde{p} = \lambda p &
\end{empheq}
\end{subequations}
with the boundary conditions
\begin{equation}
  \label{ivp-st-030}
  \begin{array}{lll}
 v(0)=p(0) &=& 0 \\
   \alpha   v_{x}(L)-\gamma  \beta  p_x(L)&=& 0 \\
  \beta  p_x(L) -\gamma  \beta  v_x(L)&=& -\frac{1}{2h^2} \lambda p(L) .
  \end{array}
\end{equation}

%First, we multiply (\ref{st-010})  by $\bar v$ and (\ref{st-020})  by $\bar p$ and then integrate by parts. Secondly, we take the conjugate of
%  both  equations (\ref{st-010}) and (\ref{st-020}) and then we multiply the conjugated equation by $v$ and  by $p$
%  respectively and integrate by parts. Then we add the results to get
%\begin{eqnarray}
%\label{star20} && \int_0^L\left( -\alpha_1  |v_x|^2 - \beta |e v_x-p_x|^2 \right)~dx - \frac{k}{h}\lambda |p(L)|^2=
%\lambda^2\int_0^L  (\rho |v|^2+ \mu |p|^2)~dx.
%\end{eqnarray}
%Now let $\lambda=\lambda_1 + i\lambda_2$ where $\lambda_1, \lambda_2\in \mathbb{R}.$ Then (\ref{star20}) has both real an imaginary components
%\begin{eqnarray}
%\label{star21}  \int_0^L\left( -\alpha_1  |v_x|^2 - \beta |e v_x-p_x|^2 \right)~dx - \frac{k}{h}\lambda_1 |p(L)|^2&=&
%(\lambda_1^2-\lambda_2^2)\int_0^L  (\rho |v|^2+ \mu |p|^2)~dx\quad\quad\quad \\
%\label{star22}   -\frac{k}{h}\lambda_2 |p(L)|^2&=& 2\lambda_1\lambda_2\int_0^L  (\rho |v|^2+ \mu |p|^2)~dx.
%\end{eqnarray}
%If $\lambda_2=0,$ then (\ref{star21}) implies that $\lambda_1<0.$ If $\lambda_2\ne 0,$ then (\ref{star22}) implies that $\lambda_1<0.$ In either case, we have that ${\rm{Re}}\lambda<0.$ Since $0\in \rho(\mc A),$ if we can show that there are no eigenvalues on the imaginary axis, or in other words, the set
Since $0\in \rho(\mc A_d ),$ if we can show that there are no eigenvalues on the imaginary axis, or in other words, the set
\begin{eqnarray}\left\{Y\in \mathrm{H} ~|~ {\rm Re} \left<\mc A_d  Y, Y\right>_{\mathrm{H}} = -\frac{1}{2h^2}|\tilde{p} (L)|^2= 0\right\}
\label{Eq81}
\end{eqnarray}
has only the trivial $Y=0$ solution, then by Arendt-Batty's stability theorem \cite{A-B},  $e^{\mc A_d (t) }$  is  a strongly stable semigroup.
Since  $\tilde{p}=\lambda p$ where $\lambda\ne 0$ by Theorem \ref{iso}, (\ref{Eq81}) implies that  $p(L)=0.$

Let $\lambda=i\tau$ where $\tau\in \mathbb{R}\backslash\{0\}.$  The eigenvalue problem (\ref{pdes-stabil10})-(\ref{ivp-st-030}) can be written
\begin{subequations}
  \label{pdes-stabil2}
\begin{empheq}[left={\phantomword[r]{0}{ }  \empheqlbrace}]{align}
\label{st-91}  &  v_{xx}=\frac{-\tau^2}{\alpha_1  }\left(\rho v + \gamma \mu p\right)   & \\
\label{st-92} &  p_{xx} =-\tau^2\left(\frac{\gamma \rho}{\alpha_1  }v+\left(\frac{\gamma^2\mu }{ \alpha_1   }+\frac{\mu}{\beta} \right)p\right) &
\end{empheq}
\end{subequations}
with the over-determined boundary conditions
\begin{eqnarray}
  \label{ivp-st-031}
 &v(0)=p(0)= p_x(L)=v_x(L)= p(L)= 0.&
\end{eqnarray}
Proving strong stability reduces to showing that (\ref{pdes-stabil2},\ref{ivp-st-031}) has only the trivial solution. Let $Z=[v, v_x, p, p_x].$ We write the system (\ref{pdes-stabil2}) in the form
\begin{eqnarray}
\label{system}\frac{d Z}{dx}=\mc D Z= \left( {\begin{array}{*{20}c}
 0 & 1 & 0 & 0  \\
 \frac{-\rho \tau^2}{\alpha_1  } & 0 & \frac{-\gamma \mu\tau^2}{\alpha_1 } & 0  \\
  0 & 0 & 0 & 1  \\
   \frac{-\gamma \rho\tau^2}{\alpha_1  } & 0 & -\left(\frac{\gamma^2\mu}{ \alpha_1  }+ \frac{\mu}{\beta}\right)\tau^2 & 0  \\
\end{array}} \right)Z.
\end{eqnarray}
The solution to (\ref{system}) is
\begin{eqnarray}
\label{system1}Z=e^{\mc D x}K
\end{eqnarray}
 where $K=[k_1, k_2, k_3, k_4]$ is the vector with arbitrary coefficients. The characteristic equation, ${\rm Det}(\mc D Z-\tilde\lambda Z)=0,$ is
$$\tilde\lambda^4 + \tilde\lambda^2\left(\frac{\gamma ^2 \mu }{\alpha_1 }+\frac{\mu}{\beta}+\frac{\rho}{\alpha_1  }\right)\tau^2 + \frac{\rho \mu}{\beta \alpha_1  }\tau^4=0.$$
This can be regarded as a quadratic equation of $\tilde \lambda^2.$ Since $\left(\frac{\gamma^2 \mu }{\alpha_1 }+\frac{\mu}{\beta}
+\frac{\rho}{\alpha_1  }\right)^2-\frac{4\rho  \mu }{\beta\alpha_1 }=\left(\frac{\gamma^2\mu }{\alpha_1 }+\frac{\mu}{\beta}
-\frac{\rho}{\alpha_1  }\right)^2+\frac{4\rho \gamma ^2\mu }{\alpha_1  ^2}>0,$ there are four roots  $\{\tilde \lambda_1^+, -\tilde \lambda_1^+, \tilde \lambda_2^{-},-\tilde \lambda_2^{-} \}$, where defining
\begin{eqnarray}
\label{lam} && a_1=\tau\zeta_1, \quad a_2=\tau\zeta_2,~~~~\tilde \lambda_1^+=ia_1,  \quad \tilde\lambda_2^-=ia_2.
\end{eqnarray}
The solution  of (\ref{system}) is written $Z=P e^{Jx} P^{-1} K$
where
\begin{eqnarray}
\label{Pmatrix} P=\left( {\begin{array}{*{20}c}
   1 & 1 & 1 & 1  \\
  i a_1 & -i a_1 & i a_2 & -i a_2  \\
      b_1 & b_1 & b_2 & b_2  \\
         i a_1 b_1 & -i a_1 b_1 & i a_2b_2 & -i a_2b_2  \\
\end{array}} \right)
\end{eqnarray}
 and $e^{Jx}={\rm{diag}}(e^{ia_1 x}, e^{-ia_1x}, e^{ia_2x}, e^{-ia_2x}),$
 \begin{eqnarray}\label{defb}b_1=\frac{1}{\gamma\mu}(\alpha_1\zeta_1^2-\rho),\quad b_2=\frac{1}{\gamma\mu}(\alpha_1\zeta_2^2-\rho),\end{eqnarray}
  or explicitly,
\begin{eqnarray}
\nonumber b_1 &=& \frac{1}{2}\left(\gamma+\frac{\alpha_1}{\gamma\beta}-\frac{\rho}{\gamma \mu}+\sqrt{\left(\gamma+\frac{\alpha_1}{\gamma\beta}-\frac{\rho}{\gamma \mu}\right)^2+\frac{4\rho}{  \mu}} \right)\\
\label{ben} b_2 &=& \frac{1}{2}\left(\gamma+\frac{\alpha_1}{\gamma\beta}-\frac{\rho}{\gamma \mu}-\sqrt{\left(\gamma+\frac{\alpha_1}{\gamma\beta}-\frac{\rho}{\gamma \mu}\right)^2+\frac{4\rho}{  \mu}} \right).
\end{eqnarray}

  Note that $b_1,b_2\ne 0,$ $b_1\ne b_2,$  and $b_1 b_2=-\frac{\rho}{\mu}.$
 The solution  of (\ref{system}) can be written   $Z=Pe^{Jx}P^{-1} K$ where
\begin{eqnarray}
\nonumber  P e^{Jx} P^{-1}=&&\left( {\begin{array}{*{20}c}
   \frac{b_1\cos{a_2 x}  -b_2\cos{a_1 x}}{b_1-b_2} & \frac{a_1 b_1\sin{a_2 x} -a_2b_2\sin{a_1 x}}{a_1 a_2(b_1-b_2)}  \\
\frac{-a_2 b_1\sin{a_2 x}  +a_1b_2\sin{a_1 x}}{b_1-b_2} & \frac{b_1 \cos{a_2 x} -b_2\cos{a_1 x}}{(b_1-b_2)} \\
   \frac{(-\cos{a_1 x}  +\cos{a_2 x})b_1b_2}{b_1-b_2} & \frac{(a_1\sin{a_2 x} -a_2\sin{a_1 x})b_1b_2}{a_1 a_2(b_1-b_2)} \\
 \frac{(a_1\sin{a_1 x}  -a_2\sin{a_2 x})b_1b_2}{b_1-b_2} & \frac{(- \cos{a_1 x} +\cos{a_2 x})b_1b_2}{(b_1-b_2)}
\end{array}} \right. \ldots \\
\nonumber  &&\quad\quad\quad\quad\quad\quad\quad\quad\left. {\begin{array}{*{20}c}
    \frac{\cos{a_1 x} -\cos{a_2 x}}{b_1-b_2} & \frac{-a_1 \sin{a_2 x} +a_2\sin{a_1 x}}{a_1 a_2(b_1-b_2)}   \\
 \frac{-a_1\sin{a_1 x} +a_2\sin{a_2 x}}{b_1-b_2} & \frac{-\cos{a_2 x} +\cos{a_1 x}}{(b_1-b_2)}\\
  \frac{b_1\cos{a_1 x} -b_2\cos{a_2 x}}{b_1-b_2} & \frac{a_2b_1 \sin{a_1 x} -a_1b_2\sin{a_2 x}}{a_1 a_2(b_1-b_2)}   \\
 \frac{-a_1b_1\sin{a_1 x} +a_2b_2\sin{a_2 x}}{b_1-b_2} & \frac{b_1\cos{a_1 x} -b_2\cos{a_2 x}}{(b_1-b_2)}
\end{array}} \right) ,
\end{eqnarray}
and $K$ is the vector of arbitrary coefficients defined in (\ref{system1}).
Note that \\${\rm{Det}} (Pe^{Jx}P^{-1})=4 a_1 a_2 (b_1-b_2)^2 \ne 0$ since $a_1, a_2\ne 0$ and $b_1\ne b_2.$  Using the boundary conditions $v(0)=p(0)$ implies  $b_1\ne b_2$ and so  $k_1=k_3=0. $
  Thus the solution of the eigenvalue problem (\ref{pdes-stabil2}) is
\begin{eqnarray}
\nonumber v(x)&=&k_2\frac{a_1 b_1\sin{a_2 x} -a_2b_2\sin{a_1 x}}{a_1 a_2(b_1-b_2)} +k_4 \frac{-a_1 \sin{a_2 x} +a_2\sin{a_1 x}}{a_1 a_2(b_1-b_2)} \\
\label{ref} p(x)&=& k_2 \frac{(a_1\sin{a_2 x} -a_2\sin{a_1 x})b_1b_2}{a_1 a_2(b_1-b_2)}  + k_4\frac{a_2b_1 \sin{a_1 x} -a_1b_2\sin{a_2 x}}{a_1 a_2(b_1-b_2)}.
\end{eqnarray}
Using the other two boundary conditions $v_x(L)=p_x(L)=0$ leads to
\begin{eqnarray}
\label{eqq1}\tilde K \left( \begin{array}{l}
 k_2 \\
 k_4 \\
 \end{array} \right)=
\left( {\begin{array}{*{20}c}
   \frac{b_1 \cos{a_2 L} -b_2\cos{a_1 L}}{(b_1-b_2)} & \frac{-\cos{a_2 L} +\cos{a_1 L}}{(b_1-b_2)}  \\
\frac{(- \cos{a_1 L} +\cos{a_2 L})b_1b_2}{(b_1-b_2)}& \frac{b_1\cos{a_1 L} -b_2\cos{a_2 L}}{(b_1-b_2)} \\
\end{array}} \right)
\left( \begin{array}{l}
 k_2 \\
 k_4 \\
 \end{array} \right)=0.\quad
\end{eqnarray}
Observe that ${\rm{Det}}\tilde K= \cos{a_1 L} \cos{a_2 L}=0$ if and only if $\cos{a_1 L} =0 $ or $\cos{a_2 L}=0. $  If for some integers $n,m,$ $a_1=\left(2n+1\right)\frac{\pi}{2L}$ and $a_2=\left(2m+1\right)\frac{\pi}{2L}$ then $v_x(L)=p_x(L)=0$ for all choices of $k_2$, $k_4$. We can choose $k_4 $ so that $p(L)=0.$ Hence the  controlled system has  an imaginary eigenvalue and is not strongly stable.   On the other hand,  if $a_1\ne \left(2n+1\right)\frac{\pi}{2L}$ and $a_2\ne\left(2m+1\right)\frac{\pi}{2L},$ then there is only the trivial solution $Z=0$. It follows from the Arendt-Batty's Theorem that the controlled system is strongly stable. Suppose now that $a_1=\left(2n+1\right)\frac{\pi}{2L}$ and $a_2\ne\left(2m+1\right)\frac{\pi}{2L}.$ Then ${\rm{Det}}\tilde K=0$ and  $k_4=b_1 k_2.$ The solution  with a parameter $k_2$ is
 \begin{eqnarray}
\nonumber v(x)&=&k_2 \frac{\sin{a_1x}}{a_1}, ~p(x)=k_2\frac{b_1\sin{a_1x}}{a_1}.
\end{eqnarray}
However, since $p(L)=0,$  the only way to obtain $p(L)=0$ is to choose $k_2=0$ and so $Z=0$.    The argument is identical if $a_1 \ne \left(2n+1\right)\frac{\pi}{2L}$ and $a_2=\left(2m+1\right)\frac{\pi}{2L}.$  Thus, if $a_1\ne\left(2n+1\right)\frac{\pi}{2L}$ or $a_2 \ne \left(2m+1\right)\frac{\pi}{2L},$ the system is strongly stable. Thus, the system is strongly stable if and only if there are not integers $n,m$ so that $\frac{a_1}{a_2} = \frac{2n+1}{2m+1} .$ This proves the theorem. $\Box$

The following theorem about the original control system (\ref{Semigroup}) is immediate.

\begin{theorem}
For any $k>0$, the control system (\ref{Semigroup}) with feedback control $V(t) = k \dot{p} (L,t) $,  is strongly stable if and only if  $\frac{\zeta_1}{\zeta_2}\ne \frac{2n-1}{2m-1},$ for some $ n,m \in \mathbb{N} . $
\end{theorem}

%*OZKAN discuss what $\dot{p} (L)$ is physically, with reference.
The  feedback signal $V(t)=k \dot p(L)=k \int_0^L  \dot D_3(\xi,  t)~ d\xi$  is physical since  $\dot p(L)$ denotes the current flowing through the electrodes of the beam  (\ref{defp}).
The limiting case of static magnetic effects corresponds to $\mu=0.$ In this case the boundary value problem (\ref{pdes-stabil2}) with over-determined boundary conditions has only the trivial solution and the controlled system is strongly stable, as is well-known.  However, in this case $B^* \varphi (x,t)=\frac{\gamma}{h}\dot v(L)$ where $\dot v(L)$ represents the velocity of the beam at  $x=L.$

Note that strong stability is  achieved with the feedback $ V(t)=k \dot{p} (L,t)$ except  for a set of coefficients $ \frac{\zeta_1}{\zeta_2}$ with Lebesgue measure zero.
% In the next section  exponential stabilizability and other values of the coefficients is considered.

%vspace{0.1in}

%%%%%%%%%
%Exact observability and exponential stabilizability
%%%%%%%%%

\section{Exact observability and exponential stabilizability} The stabilizability of the controlled piezoelectric beam will be shown to be determined by the observability of the same system.
 We start with standard definitions of exact observability, exponential stability and stabilizability, and optimizability.
\label{Sec-IV}
%vspace{0.1in}

\begin{definition} The pair  $(A, B^*)$ is exactly observable in time $T>0$  if there exists a positive constant C(T) such that for all $\varphi^0\in\mathrm H$
$$\int_0^T \|B^* e^{\mc A t} \varphi^0\|^2~ dt \ge C(T) \|\varphi^0\|^2_{\mathrm H}.$$
\end{definition}

\begin{definition}  The semigroup $\{e^{\mc A t}\}_{t\ge 0}$ with the generator $\mc A$ is exponentially stable on $\mathrm H$ if there exists constants $M, \mu>0$ such that $\|e^{\mc A t}\|_{\mathrm H}\le M e^{-\mu t}$ for all $t\ge 0.$
\end{definition}

\begin{definition} \label{a-f}  We say that the scalar $k$ is an admissible feedback  for transfer function $G(s)=C(sI-\mc A) B$
if $I-G(s) k$ has an inverse that is uniformly bounded on some right-half-plane.

\end{definition}
\begin{definition} \label{stabilizable}  The pair $(\mc A, B)$ is exponentially stabilizable on $\mathrm H$ if there exists $F\in \mc L({\rm Dom}(\mc A), \mathbb C)$ such that $(\mc A, B, F)$ is a regular triple (the transfer function $F_\Lambda(sI-\mc A)^{-1}B$ is well-defined), $1$ is an admissible feedback  for the transfer function $F_\Lambda(sI-\mc A)^{-1}B,$ and $A+B F_\Lambda$ with the domain  $\text{Dom}(A+BF_\Lambda)=\{z\in \text{Dom}(F_\Lambda)~:~ Az+BF_\Lambda z \in \mathrm{H}\}$ generates an exponentially stable semigroup $\{e^{(\mc A + B F_\Lambda) t}\}_{t\ge 0}$  on $\mathrm H.$ In the above,  the operator $F_\Lambda$ is the $\Lambda-$extension of $F:$ $$F_\Lambda z=\lim \limits_{\lambda \to \infty} F\lambda (\lambda I-A)^{-1}z$$ for all $z\in \mathrm H$ for which the limit makes sense.
\end{definition}

\begin{definition} \label{optimizable} The pair $(\mc A  , B)$ is optimizable  if
for any $z^0\in\mathrm H$ there exists a control $u(t)\in \Ltwo(0,T)$ such that  $z(t)\in \Ltwo(0,T; \mathrm H)$ where
 \begin{eqnarray}z(t)=e^{\mc A   t}z^0 + \int_0^t e^{\mc A   (t-\tau)}B u(\tau) ~d\tau. \label{gag}\end{eqnarray}
\end{definition}

%vspace{0.1in}
It is clear from the definitions that if $(\mc A , B)$ is stabilizable, then it is optimizable.   The converse of this statement  is in general false for unbounded  $B$\cite{R-W}.
%The following result is straightforward.
%\begin{lemma} \label{imply} If the pair $(\mc A , B)$ is not optimizable, then the pair $(\mc A , B)$ is not stabilizable on $\mathrm H$.
%\end{lemma}

%vspace{0.1in}

%\textbf{Proof:} By  Definitions \ref{stabilizable} and \ref{optimizable}, if the pair $(\mc A , B)$ is not optimizable, there exists a $z^0 \in \mathrm H$ such that there is no control $u(t) \in \Ltwo(0, T)$  so that $z(t) \in \Ltwo(0, T; \mathrm H)$ where $z(t)$ is defined by (\ref{gag}). This obviously implies that there is no feedback control $u(t)=\mc F z(t)$  that makes the semigroup $\{e^{(\mc A + B F) t}\}_{t\ge 0}$ exponentially stable on $\mathrm H.$ This proves the statement of the theorem.

%vspace{0.1in}
A result in \cite{Ammari-Tucsnak} implies that exact observability of the pair $(\mc A, B^*)$ in finite time on $\mathrm{H}$  is equivalent to exponential stability of the semigroup $\{e^{\mc A_d t}\}_{t\ge 0}$ on $\mathrm H$.
%vspace{0.1in}

\begin{theorem} \label{main-stab}The semigroup $\{e^{\mc A_d t}\}_{t\ge 0}$ is exponentially stable on $\mathrm H$ if and only if the pair $(\mc A, B^*)$  (defined in (\ref{Semigroup-H})) is exactly observable in finite time on $\mathrm H.$
\end{theorem}

\textbf{Proof:}
Recall  that the operator $A$ (\ref{hom-A})  is self-adjoint and positive definite and  also $B_0 \in {\mathcal L} (\mathbb C ,\mX_{-\frac{1}{2}} )$ where $\mX_1$ is  ${\rm Dom} A$ with the norm $\| A \cdot \|$  and the state-space
$\mathrm H=  \mX_{\frac{1}{2} }\times \mX. $ (\ref{duals}).Furthermore, the transfer function $G(s)$ is uniformly bounded in any right-hand plane with $\Re s \geq s_1 > 0$. Thus, the assumptions of
\cite[Thm. 2.2]{Ammari-Tucsnak} are satisfied and the conclusion follows.
$\square$
%vspace{0.1in}

%Note that the following implication
%\begin{center}``\emph{The semigroup $\{e^{\mc A_d t}\}_{t\ge 0}$ is exponentially stable in $\mathrm H$ $\Rightarrow$ the pair $(\mc A, B^*)$ is %exactly observable in $\mathrm H$ in finite time"}
% \end{center}follows from Russell's principle, i.e. stabilizability implies controllability. This does not require the condition (\ref{sok-1}). However, the other side of the implication of Theorem \ref{main-stab} (which is needed to prove Corollary {\ref{exp-stab}}) holds  true only  with the condition (\ref{sok-1}). That is why we need Lemma \ref{appen-1} for  Theorem \ref{main-stab} to hold  true.

%We also recall a result from \cite{Weiss} to establish the equivalency of optimizability of the pair $(\mc A_d , B)$ and the exponential %stability of the semigroup $\{e^{\mc A_d t}\}_{t\ge 0}.$
%%vspace{0.1in}

\begin{theorem} \label{main-optim}  The semigroup $\{e^{\mc A_d t}\}_{t\ge 0}$ is exponentially stable on $\mathrm H$ if and only if the pair $(\mc A_d , B)$ is optimizable, i.e.
for any $z^0\in\mathrm H$ there exists $u(t)\in \Ltwo(0,T)$ such that  $z(t)\in \Ltwo(0,T; \mathrm H)$ where
$z(t)$ is defined by (\ref{gag}).
\end{theorem}
%vspace{0.1in}

\textbf{Proof:} Since (\ref{Sg}) defines a well-posed and conservative system by Theorem \ref{eigens}, and it is of the class studied in \cite{Weiss}, the conclusion of the theorem follows from \cite[Thm. 1.3]{Weiss}. $\square$

%vspace{0.1in}

The following result is an immediate consequence of the two preceding results.

\begin{theorem}
\label{thm:sum}
The control system $( \mathcal A, B ) $ is optimizable if and only if $(\mathcal A, B)$ is exactly observable.
\end{theorem}

{\bf Proof:}
Theorems  \ref{main-stab} and \ref{main-optim} imply that $( \mathcal A_d, B ) $ is optimizable if and only if $(\mathcal A, B)$ is exactly observable.
Optimizability of a well-posed system is invariant under an admissible feedback  \cite[Thm. 6.3]{R-W}.  Since $(\mathcal A, B, B^*)$ is well-posed (Thm. \ref{w-pf}),
the pair $(\mc A, B) $ is optimizable if and only if $(\mathcal A, B)$ is exactly observable. $\square$

Since stabilizability implies optimizability, if  $(\mathbb A, B)$ is not exactly observable  there is no feedback controller $V(t)\in \Ltwo(0,T)$ that makes the system (\ref{Semigroup})  exponentially stable on $\mathrm{H}.$
%vspace{0.1in}
%\textbf{Notation:} Let $\mathcal{S}$ be a set, and $f,g$ be nonnegative functions on $\mathcal{S}.$ We will write $f\asymp g$ if there exists $\tilde %C_1, \tilde C_2>0$ such that
%\begin{eqnarray}\label{norm-eq}\tilde C_1 f(x) \le g(x)\le \tilde C_2 f(x), ~~ \forall x \in \mathcal{S}. \end{eqnarray}
We  now turn our attention to the observability  in the energy space $\mathrm H$ of the pair $(\mc A, B^*) .$
The following result on the  eigenvalues and eigenvectors of $\mc A$ will be  needed.

 %Let $\mathbb{Z}^*=\mathbb{Z}-\{0\}.$ The following theorem describes the eigenvalues and eigenfunctions of the operator $\mc A$ defined by (\ref{Semigroup}).

%Notation: Let $\mathcal{S}$ be a set, and $f,g$ be nonnegative functions on $\mathcal{S}.$ We will write $f\asymp g$ if there exists $C>0$ such that
%$$\frac{1}{C} f(\lambda) \le g(\lambda)\le Cf(\lambda), ~~ \forall \lambda \in \mathcal{S}. $$

%vspace{0.1in}

\begin{theorem} \label{main-thm}  Let  $\sigma_{j}=\frac{ (2j-1) \pi }{2L}, \quad j\in \mathbb{N}.$
%Assume that  $\zeta_1, \zeta_2 \in \mathbb{R}-\mathbb{Q}$ (from Theorem \ref{stronglystable}).
The  operator $\mc A$  has eigenvalues
 \begin{eqnarray}\label{EVs}\lambda^{\mp}_{1j}=\frac{\mp i\sigma_j}{\zeta_1},\quad \lambda_{2j}^{\mp}=\frac{\mp i\sigma_j}{\zeta_2}, \quad  j\in \mathbb{N}. \end{eqnarray}
The corresponding eigenfunctions are, using $\lambda_{1j}^-=-\lambda_{1j}^+, ~~\lambda_{2j}^-=-\lambda_{2j}^+,$
\begin{eqnarray}\nonumber \Psi_{1j}=\left( \begin{array}{c}
 \frac{1}{\lambda_{1j}^+} \\
 \frac{b_1}{\lambda_{1j}^+} \\
 1\\
b_1
 \end{array} \right)\sin \sigma_j x, &&~ \Psi_{-1j}=\left( \begin{array}{c}
 \frac{1}{\lambda_{1j}^+} \\
\frac{ b_1}{\lambda_{1j}^+} \\
-1 \\
-b_1
 \end{array} \right)\sin \sigma_j x,\\
\label{evectors}  ~\Psi_{2j}=\left( \begin{array}{c}
 \frac{1}{\lambda_{2j}^+} \\
\frac{ b_2}{\lambda_{2j}^+} \\
1\\
b_2
 \end{array} \right)\sin \sigma_j x, &&~\Psi_{-2j}=\left( \begin{array}{c}
 \frac{1}{\lambda_{2j}^+} \\
 \frac{b_2}{\lambda_{2j}^+} \\
-1\\
-b_2
 \end{array} \right)\sin \sigma_j x,\quad j\in \mathbb N\end{eqnarray}
 where $\zeta_1, \zeta_2, b_1, b_2$ are defined by (\ref{lam1}), (\ref{lam2}) and (\ref{defb}), respectively.
%  The corresponding eigenfunctions  $\{\sin \sigma_j x\}_{j\in\mathbb{Z}^+}$ of $\mc A$ forms an orthogonal basis in $\Ltwo(\Omega)$.

The eigenfunctions form an orthogonal basis for $\mathrm H$ and so every $\varphi^0 \in \mathrm H$ can be written, for some choice of constants $\{c_{kj}, d_{kj}\in \mathbb{C}, \quad k=1,2, \quad j\in\mathbb{N}\} ,$
 \begin{eqnarray}
 \nonumber \varphi^0&=& \sum\limits_{j \in \mathbb{N}}  \left[c_{1j}\Psi_{1j}
 + d_{1j} \Psi_{-1j}  + c_{2j}\Psi_{2j} +d_{2j}\Psi_{-2j}\right]\\
\label{init-F} &=& \sum\limits_{j \in \mathbb{N}} \left( \begin{array}{c}
\frac{1}{\lambda_{1j}^+}(c_{1j}+ d_{1j})+ \frac{1}{\lambda_{2j}^+}(c_{2j} + d_{2j}) \\
 \frac{b_1}{\lambda_{1j}^+} (c_{1j}+ d_{1j})+\frac{b_2}{\lambda_{2j}^+} (c_{2j} + d_{2j}) \\
(c_{1j}- d_{1j})+(c_{2j} -d_{2j})\\
b_1 (c_{1j}- d_{1j})+b_2(c_{2j} - d_{2j})
 \end{array} \right)\sin\sigma_j x .
\end{eqnarray}

Also,  there are  positive constants $\tilde C_1, \tilde C_2$ independent of the choice of $\Psi^0\in \mathrm H$ so that
%corresponding to the initial data (\ref{init-F}) is given
%\begin{eqnarray}\label{norm-init}\|\Psi^0\|_{\mathrm{H}}^2 \asymp \sum\limits_{j \in \mathbb{N}} \left(|c_{1j}|^2 + |d_{1j}|^2+ |c_{2j}|^2 + |d_{2j}|^2\right), ~{\rm i.e.}\end{eqnarray}
\begin{eqnarray}\label{norm-init}\tilde C_1~ \|\varphi^0\|_{\mathrm{H}}^2 \le \sum\limits_{j \in \mathbb{N}} \left(|c_{1j}|^2 + |d_{1j}|^2+ |c_{2j}|^2 + |d_{2j}|^2\right)\le \tilde C_2 ~\|\varphi^0\|_{\mathrm{H}}^2  \, . \end{eqnarray}

The function
\begin{eqnarray}
\label{sol} \varphi&=& \sum\limits_{j \in \mathbb{N}}  \left[c_{1j}\Psi_{1j}e^{\lambda_{1j}^+t}
+ d_{1j} \Psi_{-1j} e^{-\lambda_{1j}^+t} + c_{2j}\Psi_{2j} e^{\lambda_{2j}^+t}+d_{2j}\Psi_{-2j}e^{-\lambda_{2j}^+t}\right]\quad\quad\quad
\end{eqnarray}
solves   (\ref{Semigroup-H}) for the initial data (\ref{init-F}).
\end{theorem}
%By using %(\ref{Pmatrix}), we expand the initial data $(v^0, p^0, v^1, p^1)\in \mathrm{H}'$ into Fourier series
%vspace{0.1in}

\textbf{Proof:} See Appendix \ref{appenB}.
%vspace{0.1in}

We now  prove that  the pair $(\mc A, B^*)$ corresponding to (\ref{Semigroup-H}) is not exactly observable for almost all choices of parameters. The following lemma from \cite{Scott} is needed to prove this result.
%vspace{0.1in}

\begin{lemma}\label{sac} For every irrational number $\zeta$ there exists increasing sequences of coprime odd integers $\{\tilde p_m\},\{\tilde q_m\}$ and a constant $C_\zeta\ge 1$  satisfying the asymptotic relation
\begin{eqnarray}\label{irr-1}\left|~\zeta - \frac{\tilde p_m}{\tilde q_m}~\right|\le \frac{C_\zeta}{{\tilde q_m}^{2}},\quad m\to\infty.\end{eqnarray}
\end{lemma}

%vspace{0.1in}

\begin{theorem} \label{notobs}Assume that $\frac{\zeta_2}{\zeta_1}\in \mathbb{R}-\mathbb{Q}.$ Then the pair $(\mc A, B^*)$ corresponding to (\ref{Semigroup-H}) is not exactly observable on $\mathrm H.$
\end{theorem}
%vspace{0.1in}

\textbf{Proof:} Let
the sequences $\{ {\tilde p_m}\}$ and $\{ {\tilde q_m}\}$ be chosen as  in Lemma \ref{sac} with $\zeta=\frac{\zeta_2}{\zeta_1}$ and \begin{eqnarray}\label{slkk}\left|~\frac{\zeta_2}{\zeta_1}-\frac{{\tilde p_m}}{ {\tilde q_m}}\right|\le \frac{C_\zeta}{ {\tilde q_m}^{2}}.\end{eqnarray}
Define
\begin{equation}\label{kappa}\kappa_{1m}=\left\{ \begin{array}{rl}
 -  1, & ~~{\rm if} ~ {\tilde q_m}+1\equiv 0 ~({\rm {mod}}~ 4)  \\
    1, & \quad~~ {\rm otherwise}
 \end{array} \right., ~
\kappa_{2m}=\left\{ \begin{array}{rl}
 -1, &  ~~{\rm if} ~ {\tilde p_m}+1\equiv 0 ~({\rm {mod}}~ 4)  \\
  1, & \quad~~ {\rm otherwise}
 \end{array} \right.
\end{equation}
so that $\kappa_{1m} \sin \left(\frac{ {\tilde q_m} \pi}{2}\right)=\kappa_{2m} \sin \left(\frac{ {\tilde p_m} \pi}{2}\right)=1,$
and  \begin{eqnarray}\label{e-values}\lambda_{1m}=i\left({\frac{ {\tilde q_m} \pi }{2L\zeta_1}}\right), ~\lambda_{2m}=i\left({\frac{ {\tilde p_m} \pi }{2L\zeta_2}}\right).
\end{eqnarray}
Defining  \begin{eqnarray} \label{count-inits} \Phi_{1m}^0&=\frac{\kappa_{1m}}{b_1} \left( \begin{array}{c}
 \frac{1}{\lambda_{1m}}\\
\frac{b_1}{\lambda_{1m}}\\
1\\
  b_1
 \end{array} \right)\sin\left({\frac{  {\tilde q_m} \pi x}{2L}}\right),~~\Phi_{2m}^0= \frac{\kappa_{2m}}{b_2}\left( \begin{array}{c}
 \frac{1}{\lambda_{2m}}\\
\frac{b_2}{\lambda_{2m}}\\
1\\
b_2
 \end{array} \right)\sin\left({\frac{  {\tilde p_m} \pi x}{2L}}\right)\quad\quad
 \end{eqnarray}
where $b_1$ and $b_2$ are defined by (\ref{defb}),
\begin{eqnarray} \nonumber \Phi_{1m}&=&=\frac{\kappa_{1m}}{b_1} \left( \begin{array}{c}
 \frac{1}{\lambda_{1m}}\\
\frac{b_1}{\lambda_{1m}}\\
1\\
  b_1
 \end{array} \right)\sin\left({\frac{  {\tilde q_m} \pi x}{2L}}\right)e^{\lambda_{1j}t},\\
\label{count-sols} \Phi_{2m}&=& \frac{\kappa_{2m}}{b_2}\left( \begin{array}{c}
 \frac{1}{\lambda_{2m}}\\
\frac{b_2}{\lambda_{2m}}\\
1\\
b_2
 \end{array} \right)\sin\left({\frac{  {\tilde p_m} \pi x}{2L}}\right)e^{\lambda_{2j}t}
 \end{eqnarray}
are the solutions of (\ref{Semigroup-H}) with the initial conditions $\Phi_{1m}(x,0)=\Phi_{1m}^0$ and $\Phi_{2m}(x,0)=\Phi_{2m}^0,$ respectively.
(This  follows easily from (\ref{init-F}) with the choices of $c_{2j}=d_{2j}=d_{1j}\equiv 0$ for all $j\in \mathbb{N},$  $c_{1j} \equiv 0$ for $j\ne m,$ and $c_{1m}=\frac{\kappa_{1m}}{b_1}$ for the first solution,  and $c_{1j}=d_{1j}=d_{2j}\equiv 0$ for all $j\in \mathbb{N},$  $c_{2j} \equiv 0$ for $j\ne m,$ and $c_{2m}=\frac{\kappa_{2m}}{b_2}$ for the second solution.) By  linearity,   $\Phi_m=\Phi_{1m}-\Phi_{2m}$ is the solution of (\ref{Semigroup-H}) corresponding to the initial condition $\Phi_m(x,0)=\Phi_{1m}(x,0)-\Phi_{2m}(x,0)=\Phi_{1m}^0-\Phi_{2m}^0.$ Using (\ref{inner}) and (\ref{count-inits})
\begin{eqnarray}\nonumber \|\Phi_m(x,0)\|^2_{\mathrm{H}}&=&\|\Phi_{1m}^0\|_{\mathrm H}^2 + \|\Phi_{2m}^0\|_{\mathrm H}^2 \\
\nonumber &=& \frac{L}{2}\left[\frac{\rho}{b_1^2}+\frac{\rho}{b_2^2}+2\mu+\zeta_1^2\left(\frac{\alpha_1}{b_1^2} + \frac{\beta}{b_1^2}(\gamma-b_1^2)\right)+\zeta_2^2\left(\frac{\alpha_1}{b_2^2} + \frac{\beta}{b_2^2}(\gamma-b_2^2)\right)\right]\\
\label{sabit} &=& {\rm constant}.
%\nonumber &=& L\left(2\mu + \rho \left(\frac{1}{b_1^2}+\frac{1}{b_2^2}\right)\right)\\
%\nonumber  &\asymp& O(1)
% \quad {\rm as}\quad   {\tilde q_m}\to \infty.
\end{eqnarray}
 Recalling   the definition of the operator $B$ (\ref{defb_0}), (\ref{slkk}), and (\ref{kappa}), leads to
 \begin{eqnarray}
\nonumber \left| B^* \Phi_m\right| &=& \frac{1}{h}\left|~\kappa_{1m} \sin \left(\frac{ {\tilde q_m} \pi}{2}\right)e^{i\left({\frac{ {\tilde q_m} \pi t}{2L\zeta_1}}\right)}-\kappa_{2m} \sin \left(\frac{ {\tilde p_m} \pi}{2}\right)e^{i\left({\frac{ {\tilde p_m} \pi t}{2L\zeta_2}}\right)}~\right|\\
\nonumber &=&\frac{1}{h}\left|~ e^{i\left({\frac{ {\tilde q_m} \pi t}{2L\zeta_1}}\right)}- e^{i\left({\frac{ {\tilde p_m} \pi t}{2L\zeta_2}}\right)}~\right|\\
%\nonumber &=& \sqrt{~2-2\cos{\left({\frac{ {\tilde q_m} \pi t}{2L\zeta_1}}-{\frac{ {\tilde p_m} \pi t}{2L\zeta_2}}\right)}}\\
%\nonumber & =& 2\left|~\sin{\left({\frac{ {\tilde q_m} \pi t}{4L\zeta_1}}-{\frac{ {\tilde p_m} \pi t}{4L\zeta_2}}\right)}\right|\\
%\nonumber &\le& \frac{\pi t}{2L} \left|\frac{ {\tilde q_m} }{\zeta_1}-{\frac{ {\tilde p_m} }{\zeta_2}}\right|\\
 \nonumber &\le& \frac{\pi t}{2hL} \left|\frac{ {\tilde q_m} }{\zeta_1}-{\frac{ {\tilde p_m} }{\zeta_2}}~\right|\\
\label{slk}  &\le&   \frac{\pi t C_{\zeta}}{2Lh\zeta_2  {\tilde q_m}}
%\nonumber &=& \frac{\pi}{2L}\left|\left(\frac{\tilde q_j}{\zeta_1}-\frac{\tilde p_j}{\zeta_2}\right)e^{i\left({\frac{\tilde q_j \pi t}{2L\zeta_1}}\right)}-\frac{\tilde p_j}{\zeta_2}\left(e^{i\left({\frac{\tilde p_j \pi t}{2L\zeta_2}}\right)}-e^{i\left({\frac{\tilde q_j \pi t}{2L\zeta_1}}\right)}\right)\right|\\
\end{eqnarray}
where  the Mean Value Theorem  was used to obtain the third line from the second line.
Therefore, defining $M=\frac{\pi^2 T^3 C_{\zeta}^2}{12 L^2 h^2\zeta_2^2 },$
\begin{eqnarray}
\nonumber \int_0^T | B^* \Phi_m|^2 ~dt &\le & \frac{M}{ {\tilde q_m}^2} .
\end{eqnarray}
Now  $\|\Phi_m(x,0)\|_{\mathrm{H}}^2$ is constant while $\int_0^T | B^* \Phi_m|^2 ~dt =O( {\tilde q_m}^{-2}).$ Therefore the pair $(\mc A, B^*)$  is not exactly observable on $\mathrm{H}$ if $\frac{\zeta_2}{\zeta_1}\in \mathbb{R}-\mathbb{Q}.$  $\square$
%\begin{corollary}\label{rat} Let $~\frac{\zeta_2}{\zeta_1}=\frac{p}{q}$ be a rational number, then the system (\ref{homo-vol-H})  is not exactly observable.
%\end{corollary}

%\textbf{Proof:} The term
%$\int_0^T \left|\frac{ q \pi t}{2L\zeta_1}-{\frac{ p \pi t}{2L\zeta_2}}\right|^2~dt$
%in (\ref{eq31}) is identically equal to zero.
%%vspace{0.1in}

%vspace{0.1in}
\begin{corollary} \label{irr-not-stab} If $\frac{\zeta_2}{\zeta_1}\in \mathbb{R}-\mathbb{Q}$ then  $(\mc A, \mc B)$ is not exponentially stabilizable on $ \mathrm{H}.$
\end{corollary}
%vspace{0.1in}

\textbf{Proof:} Since the pair $(\mc A, B^*)$  is not exactly observable on $\mathrm H$ by Theorem \ref{notobs},  Theorem \ref{thm:sum} implies that $(\mc A, B) $ is not optimizable. Finally, since stabilizability implies optimizability,  there is no admissible feedback operator that makes the system   exponentially stable on $\mathrm{H}.$  $\square$

%vspace{0.1in}
 \begin{corollary} Let $\frac{\zeta_2}{\zeta_1}\in \mathbb{Q}$ such that $\frac{\zeta_2}{\zeta_1}=\frac{\tilde p}{\tilde q}$ where ${\rm gcd}(\tilde p,\tilde q)=1$ and $\tilde p, \tilde q$ are both odd integers. Then the pair $(\mc A, B^*)$ corresponding to (\ref{Semigroup-H}) is not exactly observable on $\mathrm{H}.$ Therefore  the system   is not exponentially stabilizable on $ \mathrm{H}.$
\end{corollary}
%vspace{0.1in}

 \textbf{Proof:}  We can  choose  $m \in \mathbb{N}$ such that  $ {\tilde q_m}=\tilde q$ and $ {\tilde p_m}=\tilde p.$ For this particular choice of $\tilde p_m$ and $\tilde q_m,$
\begin{eqnarray}\label{angut}\left|\lambda_{1m}-\lambda_{2m}\right|=\frac{\pi}{2L}\left|\frac{i \tilde q}{\zeta_1}-\frac{i \tilde p}{\zeta_2}\right|\equiv 0.
 \end{eqnarray}
 This implies that some eigenvalues coincide, and therefore there is no gap between the eigenvalues. With the identical choice of $\Phi_m^0$   (\ref{count-inits}) with $ {\tilde q_m}=\tilde q$ and $ {\tilde p_m}=\tilde p~,$ as in the proof of Thm. \ref{notobs},    $\left| B^* \Phi_m\right|\equiv 0$ by (\ref{slk}) and (\ref{angut}) while $\|\Phi_m^0\|_{\mathrm{H}}={\rm constant}$ by (\ref{sabit}). Thus $(\mc A, B^*)$ is not exactly observable on $\mathrm{H}.$ The conclusion then follows from Theorem \ref{thm:sum}  $\square$
%vspace{0.1in}

Note that if $\frac{\zeta_2}{\zeta_1}$  can be written as a ratio of odd integers, the system is not even approximately observable.

The only remaining case  to consider is  when $\frac{\zeta_2}{\zeta_1}$ can be written as a ratio of coprime integers where one is odd and one is even. In this case eigenvalues (\ref{EVs}) have a uniform gap and the system is exactly observable.
We will use the following theorem.

%vspace{0.1in}

\begin{theorem} (Ingham's Theorem) \label{Ingham} \cite[page 162]{Young} If the strictly increasing sequence $\{s_n\}_{n\in \mathbb{N}}$ of real numbers satisfies the gap condition
\begin{eqnarray}s_{n+1}-s_n\ge \gamma\label{gap}\end{eqnarray}
for all $n \in \mathbb{N},$ for some $\gamma > 0,$ then there exists positive constants $\tilde c_3(T)$ and $\tilde c_4(T)$ such that for all $T > \frac{2\pi}{\gamma}$
  \begin{eqnarray}\label{obs-ing} \tilde c_3(T) \sum_{n\in \mathbb{N}}|g_n|^2 \le \int_0^T \left|\sum_{n\in \mathbb{N}}g_n e^{is_n t}\right|^2~dt \le \tilde c_4(T) \sum_{n\in \mathbb{N}}|g_n|^2\end{eqnarray}
\end{theorem}
for all functions $\sum\limits_{n\in \mathbb{N}}  g_n e^{is_n t}~:~ \sum\limits_{n \in \mathbb{N}} {\left| {g_n } \right|^2 }<\infty. $
%Note that the set $\tilde{\mathbb{Q}}$ not only contains the set of rational numbers $\mathbb{Q},$ but also the Liouville's numbers chosen as in %Definition \ref{irr1}.

%vspace{0.1in}

 \begin{theorem} \label{main-obs}Let $\frac{\zeta_2}{\zeta_1}\in \mathbb{Q}$ such that $\frac{\zeta_2}{\zeta_1}=\frac{\tilde p}{\tilde q}$ where ${\rm gcd}(\tilde p, \tilde q)=1$ and $\tilde p, \tilde q$ are even and odd integers, respectively; or the other way around. Choose $$T > 2L~{\min{\left(\zeta_1,\zeta_2, 2\tilde q \zeta_2 \right)}} .$$ Then the pair $(\mc A, B^*)$ is exactly observable on $\mathrm H,$ i.e. there exists a constant $C(T)>0$ such that solutions $\varphi$ of the system (\ref{Semigroup-H}) satisfy the following observability estimate:
\begin{eqnarray}
\label{obs10} & \int_0^T |B^*\varphi|^2 ~dt \ge  C(T) \|\varphi^0\|_{\mathrm{H}}^2.&
\end{eqnarray}
 \end{theorem}

\textbf{Proof:} Let $s_{1j}=\frac{\sigma_j}{\zeta_1}=\frac{(2j-1)\pi}{2L\zeta_1}$ and $s_{2j}=\frac{\sigma_j}{\zeta_2}=\frac{(2j-1)\pi}{2L\zeta_2}$ for $j\in\mathbb{N}.$ The set of eigenvalues (\ref{EVs}) can be rewritten as
\begin{eqnarray}\label{forgap}  \lambda_{kj}^{\mp}=\mp is_{kj}, ~~k=1,2, ~~ j\in \mathbb{N}.
%\nonumber &=&i\{\mp\frac{(2j-1)\pi}{2L\zeta_1}, \mp\frac{(2j-1)\pi}{2L\zeta_2}~~k=1,2, ~~ j\in \mathbb{N} \}.
\end{eqnarray}
Letting  $\varphi$ be any solution to (\ref{Semigroup-H}),  with initial condition expanded as in   (\ref{defb_0}). By  (\ref{EVs})-(\ref{norm-init})
$$
|B^* \varphi| = \left|\frac{1}{h}\sum\limits_{j \in \mathbb{N}}  \left[b_1\left(c_{1j}e^{is_{1j}t}
 - d_{1j} e^{-is_{1j}t}\right) + b_2\left(c_{2j} e^{is_{2j}t}-d_{2j}  e^{-is_{2j}t}\right)\right] (-1)^j \right|.\quad\quad\quad
$$
Showing (\ref{obs10}) is equivalent to finding a constant $C(T)$, independent of the initial condition,  so that
\begin{eqnarray}
\nonumber  && \frac{1}{h^2}\int_0^T \left|\sum\limits_{j \in \mathbb{N}}  \left[b_1\left(c_{1j}e^{is_{1j}t}
 - d_{1j} e^{-is_{1j}t}\right) + b_2\left(c_{2j} e^{is_{2j}t}-d_{2j}  e^{-is_{2j}t}\right)\right] (-1)^j \right|^2~dt\quad\\
\label{suck22} && \quad\quad\quad \ge C(T) \|\varphi^0\|_{\mathrm{H}}^2.
\end{eqnarray}

We first show that the gap condition (\ref{gap}) in Theorem \ref{Ingham} holds. If $k=n,$
\begin{eqnarray}|{\rm }s_{kj}-s_{nm}|\ge\left\{
                                   \begin{array}{ll}
                                     \frac{\pi}{L\zeta_1}, & \quad\hbox{$k=n=1$;} \\
                                    \frac{\pi}{L\zeta_2}, & \quad\hbox{$k=n=2$.}
                                   \end{array}
                                 \right., \quad j,m\in\mathbb{N}
\label{est1}
\end{eqnarray}
by (\ref{forgap}). Now let $k\ne n.$ Without loss of generality,  assume that $\tilde p$ is even and $\tilde q$ is odd. By (\ref{forgap})
\begin{eqnarray}\nonumber |s_{1j}-s_{2m}|&=&\frac{\pi}{2L}\left|\frac{2j-1}{\zeta_1}-\frac{2m-1}{\zeta_2}\right|\\
\nonumber &=&\frac{\pi}{2L}\frac{1}{\zeta_2 \tilde q}\left|(2j-1)\tilde p -(2m-1)\tilde q\right|\\
\label{est2} &\ge& \frac{\pi}{2L}\frac{1}{\zeta_2 \tilde q}\quad\quad\quad\quad\end{eqnarray}
using $|(2j-1)\tilde p -(2m-1)\tilde q|\ge 1$ since $(2j-1)\tilde p$ is an even number and $(2m-1)\tilde q$ is an odd number. Similarly $|s_{2j}-s_{1m}|\ge \frac{\pi}{2L}\frac{1}{\zeta_2 \tilde q}.$

Let's rearrange the set $\{\mp s_{kj}:~ k=1,2, ~ j\in \mathbb{N}\}$ into an increasing sequence of $\{s_n, ~n\in \mathbb{N}\},$ and denote the coefficients $\{ (-1)^j b_k c_{kj}, (-1)^{j+1}b_k d_{kj}~ :~ k=1,2, ~ j\in \mathbb{N}\}$ by $\{g_n, ~n\in\mathbb{N}\}.$   Then (\ref{est1}) and (\ref{est2}) yields $$s_{n+1}-s_{n}\ge \gamma := \frac{\pi}{L}{\min{\left(\frac{1}{\zeta_1},\frac{1}{\zeta_2}, \frac{1}{2 \zeta_2\tilde q}\right)}}, $$
and therefore  the gap condition (\ref{gap}) holds. By Theorem \ref{Ingham}, for $T>\frac{2\pi}{\gamma}=2L~{\min{\left(\zeta_1,\zeta_2, 2 \zeta_2\tilde q\right)}}$
\begin{eqnarray}
\nonumber \int_0^T |B^*\varphi|^2 ~dt &=& \frac{1}{h^2}\int_0^T \left| \sum\limits_{j \in \mathbb{N}}  \left[b_1\left(c_{1j}e^{is_{1j}t}
 - d_{1j} e^{-is_{1j}t}\right) \right.\right. \\
 \nonumber &&\left.\left. \quad\quad\quad\quad\quad\quad+~ b_2\left(c_{2j} e^{is_{2j}t}-d_{2j}  e^{-is_{2j}t}\right)\right] (-1)^j \right|^2dt\\
%\nonumber &\ge& \sum\limits_{j =  1 }^\infty \frac{1}{|\lambda_j|^2}\left(b_1^2 \left[c_{1j}^2 + d_{1j}^2\right] + b_2^2 \left[c_{2j}^2 + d_{2j}^2\right] \right)\\
\nonumber &=& \frac{1}{h^2} \int_0^T \left|~\sum\limits_{j\in \mathbb{N}} g_n e^{is_n t}~ \right|^2~dt\\
\nonumber & \ge & \tilde c_3(T)\sum\limits_{n\in \mathbb{N}}\left|g_n\right|^2\\
\nonumber &=& \tilde c_3(T)\sum\limits_{j\in \mathbb{N}}\left(b_1^2\left(\left|c_{1j}\right|^2 + \left|d_{1j}\right|^2\right) + b_2^2\left(\left|c_{2j}\right|^2 + \left|d_{2j}\right|^2\right)\right) \\
\nonumber&\ge & \tilde c_3(T) \min (b_1^2, b_2^2)\sum\limits_{j \in \mathbb{N}} \left(|c_{1j}|^2 + |d_{1j}|^2+ |c_{2j}|^2 + |d_{2j}|^2\right)\\
\label{suck21}  &\ge& C(T) \|\varphi^0\|^2_{\mathrm H}
\end{eqnarray}
where $C(T)=\frac{\tilde c_3(T)\min (b_1^2, b_2^2)}{\tilde C_1}.$ The constants $\tilde C_1$ and ${\tilde c}_3(T)$  are due to (\ref{norm-init}) and (\ref{obs-ing}), respectively. Hence (\ref{suck22}) holds and the system is exactly observable. $\square$
%vspace{0.1in}

\begin{corollary} \label{exp-stab} Let $\frac{\zeta_2}{\zeta_1}\in \mathbb{Q}$ such that $\frac{\zeta_2}{\zeta_1}=\frac{\tilde p}{\tilde q}$ where ${\rm gcd}(\tilde p, \tilde q)=1$ and $\tilde p, \tilde q$ are even and odd integers, respectively; or the other way around. The semigroup $\{e^{\mc A_d t}\}_{t\ge 0}$ is exponentially stable on $ \mathrm{H}$, and so $(\mc A , B )$ is exponentially stabilizable on $\mathrm H$.
\end{corollary}
%vspace{0.1in}

\textbf{Proof:} The pair $(\mc A, B^*)$ is exactly observable on $\mathrm{H}$ by Theorem \ref{main-obs}. Therefore the semigroup $\{e^{\mc A_d t}\}_{t\ge 0}$ is exponentially stable on $\mathrm{H}$  by Theorem \ref{main-stab}. $\square$
%vspace{0.1in}

\section{Conclusions}

The main result of this paper is to show that magnetic effects in piezoelectric beams, even though small, have a dramatic effect on observability and stabilizability. The piezoelectric beam model, without magnetic effects is exactly observable and exponentially stabilizable, by $-B^*$.  However, when magnetic effects are included, the beam is only observable and stabilizable when the parameter $\frac{\zeta_2}{\zeta_1}$ is coprime ratio of odd and even integers. In this case, the beam can be stabilized by the feedback $-B^*$. If this parameter is an  irrational number, the beam can be strongly stabilized by the $-B^*$ feedback, but not exponentially stabilized. Explicit polynomial estimates for this situation have been obtained  \cite{ozkan} .
%depends on the degree of the irrationality of the parameter $\frac{\zeta_2}{\zeta_1}$.
%The explicit polynomial decay estimates, which are valid only for more regular initial data,  are given by using Diophantine's approximations and the H\"{o}lder's inequality for weighted spaces. This and the other cases are summarized in Table \ref{results}.
%The cases where the exponential stabilizability fails but strong stability holds for $-B^*$ feedback are studied in \cite{ozkan}.

Another difference between the model with magnetic effects and without is the physical nature of the  $-B^*$ feedback.  In  models without magnetic effects (\ref{or})  this observation corresponds to the measurement of velocity of the beam at the end. However,  in the model with magnetic effects (\ref{homo-vol}) $B^*$ corresponds to the total current at the electrodes. It is typically easy to measure the current at the electrodes, much easier  than to measure velocity.

However, voltage-controlled systems exhibit hysteresis when they are actuated at high-frequencies \cite[e.g.]{Moheimani}.
Experimental evidence shows that current and charge actuation leads to much less hysteresis than voltage actuation, for instance see \cite{F-M, Furutani,Newcomb-Flinn}. A model for piezo-electric beams with magnetic effects and  current control has been derived \cite{O-M}.
The model for current control is quite different and the control operator is bounded.

No damping was considered in this paper. Including damping would of course make the system stable.
% If we add a viscous damping to the mechanical equation, that is $k\dot v$ with  $k>0,$ it can be easily proved that both systems with and without magnetic effects, (\ref{homo-vol}) and (\ref{or}) respectively, are exponentially stable even without the boundary feedback.
However, the electrical  nature, as opposed to mechanical, of $B^*$ would still remain, as would the basic conclusions of  the restricted effectiveness  of control.  As noted at the end of Section \ref{Sec-II}, modifying the Euler-Bernoulli beam to a Mindlin-Timoshenko beam makes no fundamental difference to the model since the bending and rotation parts of the model are decoupled from the stretching.

The extension to including magnetic effects in structures with piezoelectric patches is studied in \cite{accpaper} for both Euler-Bernoulli and Mindlin-Timoshenko beam models. For patches, bending and rotation equations are  coupled to the stretching equation.
%Without the magnetic effects, the regularity  and  controllability problems for only the bending equation were studied for the beam model \cite{Tucsnak-a}, and for the plate model where the patches are located on a curve on the plate (so the voltage  control is acting on a curve)
Previous research on control of structures with piezo-electric patches, without magnetic effects,  \cite{J-T,Tucsnak-a,Tucsnak-b}
 showed that the location of the patch(es) on the beam/plate strongly determines the controllability/stabilizability. The recent research discussed in \cite{accpaper} and \cite{J-T,Tucsnak-a,Tucsnak-b} suggest that  controllability/stabilizability depends on not only the location of the patches but also the system parameters.
This  is currently being studied.
%\kchg{Is the above citing correct? Or do they not all mention the importance of patch location? What about size?}
%\chg{Size of the patch(es) has not mentioned in these references. The discussion centers around where the edges of the patch(es) should be located relying on number theoretic results, once the length of the patch(es) is specified.
%Also I am not sure how the model with a patch located on a curve is derived and how physical it is. In fact, how will the piezo-patches stretch or bend in that case?}

\begin{table}[h]
\renewcommand{\arraystretch}{1.3}
\begin{center}
\begin{tabular}{|p{3cm}|c|c|c|}
  \hline
   Parameter $\frac{\zeta_2}{\zeta_1}$ & Strongly Stabzble & Exactly Obs. & Exp. Stabzble\\
\hline
  irrational & \cmark & X & X \\
\hline
$\frac{\tilde p}{\tilde q}, ~~\tilde p, \tilde q$ odd  & X & X & X  \\
\hline

$\frac{\tilde p}{ \tilde q}, ~~\tilde p$ odd , $\tilde q$ even   & \cmark & \cmark & \cmark\\
\hline
$\frac{\tilde p}{\tilde q}, ~~\tilde p$ odd , $\tilde q$ even   & \cmark & \cmark & \cmark\\
  \hline
\end{tabular}
\end{center}
\caption{Summary of results. In the first column,  $(\tilde{p}, \tilde{q})$  are coprime integers, except for the first line where  $\frac{\zeta_2}{\zeta_1}$  is irrational. }
\label{results}
%\centering \small
\end{table}

\section{Acknowledgement}
This  research was supported by a Discovery Grant from the Natural Sciences and Engineering Research Council of Canada (NSERC).

\appendix
%vspace{0.1in}

\section{Proof of Lemma \ref{appen-1}}
\label{appenB}
Our goal is to prove that  for any real $\tilde{s} >0$,    there is a constant $M>0$ such that   for all $\Re s \geq \tilde{s}, $
\begin{eqnarray}\label{sonuc} | \mb G( s ) | = |s B_0^* ( s ^2 I + A)^{-1}B_0 | \le M .  \end{eqnarray}

The transfer function $G(s)$ can be found as the  solution to the elliptic problem corresponding to the boundary control problem (\ref{homo-vol}); see, for instance, \cite{C-M}.
Define for any scalar $V $
$$\left( \begin{array}{c}
 Y \\
 Z \\
 \end{array} \right)=( s ^2 I + A)^{-1}B_0 V$$ where $Y$ and $Z$ satisfy
 %so that by (\ref{dom-hom-A}) and (\ref{defb_0})
 %\begin{eqnarray}\label{aux}( s ^2 I + A)
%\left( \begin{array}{l}
% Y \\
% Z \\
% \end{array} \right)=B_0 V=\left( \begin{array}{c}
% 0 \\
% \delta(x-L) \\
% \end{array} \right) V
% \end{eqnarray}
% with the boundary conditions
%\begin{eqnarray}
%\label{aux-b-1}
% &Y(0)=Z(0)= \alpha   Y_{x}(L)-\gamma\beta  Z_x(L)=\beta  Z_x(L) -\gamma\beta Y_x(L)= 0.&
%\end{eqnarray}
 %Solving the system (\ref{aux}) is equivalent to solving  the following system
\begin{subequations}
\label{son}
\begin{empheq}[left={\phantomword[r]{0}{ }  \empheqlbrace}]{align}
 &\alpha   Y_{xx}-\gamma \beta  Z_{xx} = \rho s ^2 Y   & \\
 &   \beta  Z_{xx} -\gamma \beta Y_{xx}= \mu  s ^2 Z , &
\end{empheq}
\end{subequations}
with the boundary conditions
\begin{eqnarray}
\label{aux-b}
 &Y(0)=Z(0)= \alpha   Y_{x}(L)-\gamma\beta  Z_x(L)=\beta  Z_x(L) -\gamma\beta Y_x(L)+ \frac{V}{h}= 0.&
\end{eqnarray}
The system (\ref{son})-(\ref{aux-b}) is similar to the system (\ref{pdes-stabil10})-(\ref{ivp-st-030})
with a slight change in the boundary conditions (\ref{ivp-st-030}). We follow the same approach  used in the proof of Theorem \ref{stronglystable} to solve the system (\ref{son}). By (\ref{ref}),  the first two boundary conditions in (\ref{aux-b}) yield the general solution
\begin{eqnarray}
\nonumber Y&=&k_1\frac{\zeta_1 b_1\sinh{(   \zeta_2 s x)} -\zeta_2 b_2\sinh{(   \zeta_1 s x)}}{  \zeta_1 \zeta_2(b_1-b_2)} +k_2 \frac{- \zeta_1 \sinh{( \zeta_2 s x)} + \zeta_2\sinh{(   \zeta_1 s x)}}{   \zeta_1\zeta_2(b_1-b_2)} \\
\nonumber Z&=& k_1 \frac{( \zeta_1\sinh{(   \zeta_2 s x)} - \zeta_2\sinh{(   \zeta_1 s x)})b_1b_2}{   \zeta_1\zeta_2(b_1-b_2)}  + k_2\frac{ \zeta_2 b_1 \sinh{(   \zeta_1 s x)} -\zeta_1 b_2\sinh{( \zeta_2  sx)}}{   \zeta_1\zeta_2(b_1-b_2)}\quad\quad\quad
\end{eqnarray}
with two arbitrary constants $k_1$ and $k_2.$ In the above $\zeta_1, \zeta_2, $  $b_1,$ and $b_2$ are the same nonzero constants defined by (\ref{lam1}), (\ref{lam2}), and (\ref{defb}), respectively. $k_1$ and $k_2$ are determined by applying the last two boundary conditions in (\ref{aux-b})
\begin{eqnarray}
\nonumber k_1 &=&\frac{\frac{  V}{\alpha_1 s }}{\cosh{\zeta_1 L}\cosh{\zeta_2 L}}\frac{
 \gamma\left(b_2 \cosh (\zeta_2 L)-b_1 \cosh (\zeta_1 L)\right) + \frac{\alpha}{\beta}\left( \cosh (\zeta_1 L)- \cosh (\zeta_2 L)\right)}{b_1-b_2}\\
 \nonumber k_2 &=&\frac{\frac{  V}{\alpha_1 s}}{\cosh{\zeta_1 L}\cosh{\zeta_2 L}}\frac{
\frac{\alpha}{\beta} \left(b_2 \cosh (\zeta_1 L)-b_1 \cosh (\zeta_2 L)\right) + \gamma\left( \cosh (\zeta_2 L)- \cosh (\zeta_1 L)\right)}{b_1-b_2}.
\end{eqnarray}
After simplifications,
$$Z(x)=\frac{V}{\alpha_1  h(b_1-b_2) s }\left[\frac{b_2(\frac{\alpha}{\beta}-b_1 \gamma)}{\zeta_2}\frac{\sinh{(  \zeta_2 s x)}}{\cosh{(  \zeta_2 s L )}}+\frac{b_1(b_2 \gamma-\frac{\alpha}{\beta})}{\zeta_1}\frac{\sinh{(   \zeta_1 s L )}}{\cosh {(   \zeta_1 s x)}}\right]$$
and therefore
\begin{eqnarray}
\nonumber \mb G( s )&=& \frac{s}{V}B_0^* Z\\
 \nonumber &=&  -\frac{s}{V h}  Z(L)\\
\label{aux-c} &=& \frac{1}{\alpha_1 h^2(b_1-b_2) }\left[\frac{b_2(b_1 \gamma-\frac{\alpha}{\beta})}{\zeta_2}\tanh{(  \zeta_2 s L)}-\frac{b_1(b_2 \gamma-\frac{\alpha}{\beta})}{\zeta_1}\tanh{(   \zeta_1 s L)}\right].\quad\quad
\end{eqnarray}

Now   bounds for the functions $\left|\tanh{ s  \zeta_1 L}\right| $ and  $\left|\tanh{ s  \zeta_2 L}\right|$ are calculated. Writing $ s = s _1+i s _2$ where $s_1 \geq \tilde{s}$ for some real $\tilde{s}>0$,
%  After simplifications we obtain
%\begin{eqnarray}\nonumber |\tanh( s  \zeta_1 L)|&=& |\tanh( s _1 \zeta_1 L + i s _2 \zeta_1 L)|\\
%\nonumber &=& \left|\frac{e^{4 s _1\zeta_1L}-1+2 i e^{2 s _1\zeta_1L}\sin (2 s _2\zeta_1L) }{e^{4 s _1\zeta_1L}+1+2 e^{2 s _1\zeta_1L}\cos (2 s _2\zeta_1L) }\right| \\
%\nonumber &=& \frac{\sqrt{\left(e^{4 s _1\zeta_1L}-1\right)^2 + 4 e^{4 s _1\zeta_1L}\sin^2 (2 s _2\zeta_1L)}}{\left|e^{4 s _1\zeta_1L}+1+2 e^{2 s _1\zeta_1L}\cos (2 s _2\zeta_1L)\right|}\\
%\nonumber &\le& \frac{e^{4 s _1\zeta_1L}+1} {\left(e^{2 s _1\zeta_1L}-1\right)^2}\\
%\nonumber &=& \frac{1}{2}\cosh  (2 s _1\zeta_1 L) ~{\rm cosech}^2 ( s _1\zeta_1L)\\
%\nonumber &:=& M_1( s _1).
%W\end{eqnarray}
% Similarly,
%$|\tanh( s  \zeta_2 L)|\le M_2( s _1):=\frac{e^{4 s _1\zeta_2L}+1} {\left(e^{2 s _1\zeta_2L}-1\right)^2}.$
%The functions $M_1~:~ \mathbb{R}^+ \to \mathbb{R}^+$ and $M_2~:~ \mathbb{R}^+ \to \mathbb{R}^+$ are continuous, strictly decreasing and convex functions. For  any $\tilde s>0$ we have
%   $$|M_1|\le \frac{e^{4 \tilde s _1\zeta_1L}+1} {\left(e^{2 \tilde s _1\zeta_1L}-1\right)^2}<\infty,\quad |M_2|\le \frac{e^{4 \tilde s _1\zeta_2L}+1} {\left(e^{2 \tilde s _1\zeta_2L}-1\right)^2}<\infty.$%
\begin{eqnarray*}
\left|\tanh ( \zeta_1 s L ) \right|&=& \left| \frac{e^{s \zeta_1 L } - e^{ s \zeta_1 L}  }{ e^{s \zeta_1 L } + e^{ s \zeta_1 L }} \right| = \frac{ \left| 1-  e^{-2 s \zeta_1 L}  \right| }{ \left| 1+  e^{ -2s \zeta_1 L }  \right|} \leq  \frac{2}{1-e^{-2s_1 \zeta_1 L } } \leq  \frac{2}{ 1-e^{-2\tilde{s} \zeta_1 L} } .
\end{eqnarray*}
A similar bound holds for  $ \left|\tanh ( \zeta_2 s L ) \right| .$
   Finally, since $\frac{b_2(b_1 \gamma -\frac{\alpha}{\beta})}{\alpha_1h^2\zeta_2(b_1-b_2) }$ and  $\frac{b_1(b_2 \gamma -\frac{\alpha}{\beta})}{\alpha_1h^2\zeta_1(b_1-b_2) }$ in (\ref{aux-c}) are all nonzero constants, there exists a positive constant $M(\tilde s)<\infty$  such that (\ref{sonuc}) holds.

Multiplying both numerator and denominator in $\bf G(s)$ by  $\gamma$ and noting that in the case of $\gamma=0,$ i.e. the system (\ref{son})-(\ref{aux-b}) is completely decoupled, $\zeta_1=\sqrt{\frac{\rho}{\alpha_1}},$ $\zeta_2=\sqrt{\frac{\mu}{\beta}},$ $b_1\gamma=0,$ $b_2 \gamma=\frac{\alpha}{\beta}-\frac{\rho}{\mu}$ where $\alpha_1=\alpha.$ Then the transfer function (\ref{aux-c}) of the decoupled system is ${\bf G}(s)=\frac{1}{h^2\sqrt{\beta\mu}}\tanh {\left( \sqrt{\frac{\mu}{\beta}}Ls\right)},$ the same transfer function as obtained from (\ref{son})-(\ref{aux-b}) with $\gamma =0.$ $\square$
\vspace{0.1in}

{\bf{Proof of Theorem \ref{main-thm}:}}   Let $\Psi=[~z_1, z_2, z_3, z_4~]^{\rm T}.$ Solving  the eigenvalue problem $\mc A \Psi= \lambda \Psi$  corresponding to (\ref{Semigroup-H}) is equivalent to solving
$$\left( \begin{array}{l}
 z_3 \\
 z_4 \\
 \end{array} \right)=\lambda \left( \begin{array}{l}
 z_1 \\
 z_2 \\
 \end{array} \right), \quad -A\left( \begin{array}{l}
 z_1 \\
 z_2 \\
 \end{array} \right)=-\lambda \left( \begin{array}{l}
 z_3\\
 z_4 \\
 \end{array} \right) = -\lambda^2 \left( \begin{array}{l}
 z_1\\
 z_2 \\
 \end{array} \right).$$
  By using $\alpha=\alpha_1 + \gamma^2\beta,$  (\ref{coef}) and (\ref{dom-hom-A}), the eigenvalue problem can be rewritten
\begin{subequations}
  \label{pdes-stabil222}
\begin{empheq}[left={\phantomword[r]{0}{ }  \empheqlbrace}]{align}
\label{macka} & z_3=\lambda z_1,\quad z_4=\lambda z_2&\\
\label{st-911}  &  z_{1xx}=\frac{\lambda^2}{\alpha_1  }\left(\rho z_1 + \gamma \mu z_2\right)   & \\
\label{st-921} &  z_{2xx} =\lambda^2\left(\frac{\gamma \rho}{\alpha_1  }z_1+\left(\frac{\gamma ^2\mu }{ \alpha_1   }+\frac{\mu}{\beta} \right)z_2\right) &
\end{empheq}
\end{subequations}
with the  boundary conditions
\begin{eqnarray}\label{bc}z_1(0)=z_2(0)=z_{1x}(L)=z_{2x}(L)=0.\end{eqnarray}
First,  find the eigenvalues. Since the solution of (\ref{pdes-stabil222}) with $\lambda=0$ is $z_1=z_2=z_3=z_4\equiv 0$;  $\lambda=0$ is not an eigenvalue.

Define
\begin{equation}z_{1j}=f_j \sin \sigma_{j} x, \quad z_{2j}=g_j \sin \sigma_{j} x, \quad \sigma_j=\frac{(2j-1)\pi}{2L}, \quad j\in \mathbb{N}.\label{modal}\end{equation}
Solutions of this form satisfy all the homogeneous boundary conditions (\ref{bc}). We seek
$f_j, g_j$ and $\lambda_j$  so that the system (\ref{pdes-stabil222}) is satisfied.
Upon substitution of  (\ref{modal}) into (\ref{pdes-stabil222}) we obtain
\begin{subequations}
  \label{pdes-stabil3}
\begin{empheq}[left={\phantomword[r]{0}{ }  \empheqlbrace}]{align}
\nonumber  &  -\sigma_j^2 f_j =\frac{\lambda^2}{\alpha_1  }\left(\rho f_j + \gamma \mu g_j \right)   & \\
\nonumber  & -\sigma_j^2 g_j =\lambda^2\left(\frac{\gamma \rho}{\alpha_1  } f_j+\left(\frac{\gamma ^2\mu }{ \alpha_1   }+\frac{\mu}{\beta} \right)g_j\right). &
\end{empheq}
\end{subequations}
Letting  $y_j=\frac{\sigma_j^2}{\lambda^2},$  this linear system has nontrivial solutions if and only if the following characteristic equation is satisfied:
$$y_j^2 + \left(\frac{\gamma ^2 \mu }{\alpha_1 }+\frac{\mu}{\beta}+\frac{\rho}{\alpha_1  }\right)y_j + \frac{\rho \mu}{\beta \alpha_1  }=0.$$
Since $\left(\frac{\gamma ^2 \mu }{\alpha_1 }+\frac{\mu}{\beta}
+\frac{\rho}{\alpha_1  }\right)^2-\frac{4\rho  \mu }{\beta\alpha_1 }=\left(\frac{\gamma ^2\mu }{\alpha_1 }+\frac{\mu}{\beta}
-\frac{\rho}{\alpha_1  }\right)^2+\frac{4\rho \gamma ^2\mu }{\alpha_1  ^2}>0,$  the roots are
$ y_{1j}=-\zeta_1^2, ~~y_{2j}=-\zeta_2^2 $ where $\zeta_1, \zeta_2\in \mathbb{R}$ are defined by (\ref{lam1}) and (\ref{lam2}), respectively. Therefore $\lambda^{\mp}_{ 1j}= \mp \frac{i\sigma_j}{\zeta_1}$ and  $\lambda^{\mp}_{2j}=\mp\frac{i\sigma_j}{\zeta_2}$ for $j\in \mathbb{N},$  and (\ref{EVs}) follows. Observe that
\begin{eqnarray}\label{minpl}\lambda_{1j}^-=-\lambda_{1j}^+, ~~\lambda_{2j}^-=-\lambda_{2j}^+. \quad j\in\mathbb{N}.\end{eqnarray}

%Next, we find the eigenvectors (\ref{evectors}). Let $\lambda^{\mp}_{1j}=\frac{\mp i\sigma_j}{\zeta_1}, ~~j\in\mathbb{N}.$
Setting $f_j=\frac{1}{\lambda^{+}_{1j}}$  yields $g_j= \frac{b_1}{\lambda^{+}_{1j}}.$ By (\ref{macka}) and (\ref{minpl}), the first two sets of eigenvectors $\Psi_{1j}, \Psi_{-1j}$ are, using the
fact that $\lambda_{1j}^-=-\lambda_{1j}^+$
$$
\nonumber  \Psi_{1j}=\left( \begin{array}{r} \frac{1}{\lambda^{+}_{1j}}\sin \sigma_{j} x \\ \frac{b_1}{\lambda^+_{1j}}\sin \sigma_{j} x \\  \sin \sigma_{j} x \\ b_1\sin \sigma_{j} x \end{array} \right) , \quad \quad
\nonumber  \Psi_{-1j}=\left( \begin{array}{r} \frac{1}{\lambda^+_{1j}}\sin \sigma_{j}  \\ \frac{b_1}{\lambda^+_{1j}}\sin \sigma_{j} x \\ -\sin \sigma_{j} x \\ -b_1\sin \sigma_{j} x  \end{array} \right)  , ~~ j\in\mathbb{N}.
$$
%
%\begin{eqnarray}\nonumber && z_{1j}=\frac{1}{\lambda^{+}_{1j}}\sin \sigma_{j} x, ~~ z_{2j}=\frac{b_1}{\lambda^+_{1j}}\sin \sigma_{j} x, ~~ z_{3j}=\sin \sigma_{j} x,~~ z_{4j}=b_1\sin \sigma_{j} x,~~ j\in\mathbb{N},\\
%\nonumber  &&z_{1j}=\frac{1}{\lambda^+_{1j}}\sin \sigma_{j} x, ~~ z_{2j}=\frac{b_1}{\lambda^+_{1j}}\sin \sigma_{j} x, ~~ z_{3j}=-\sin \sigma_{j} x,\quad z_{4j}=-b_1\sin \sigma_{j} x,~~ j\in\mathbb{N}.
%\end{eqnarray}
% Now let and $\lambda_{2j}^{\mp}=\frac{\mp i\sigma_j}{\zeta_2},~~j\in\mathbb{N}.$
Similarly, setting $f_j=\frac{1}{\lambda_{2j}^{\mp}}$ yields $g_j=\frac{b_2}{\lambda_{2j}^{\mp}}.$ By (\ref{macka}) and (\ref{minpl}), the last two eigenvectors $\Psi_{2j}, \Psi_{-2j}$ in (\ref{evectors})are
$$
\Psi_{2j} = \left( \begin{array}{r} \frac{1}{\lambda_{2j}^+}\sin \sigma_{j} x \\ \frac{b_2}{\lambda_{2j}^+}\sin \sigma_{j} x \\ \sin \sigma_{j} x \\ b_2   \sin \sigma_{j} x \end{array} \right) , \quad \quad
\Psi_{-2j} = \left( \begin{array}{r}  \frac{1}{\lambda_{2j}^+}\sin \sigma_{j} x \\ \frac{b_2}{\lambda_{2j}^+}\sin \sigma_{j} x \\ -\sin \sigma_{j} x\\ -b_2   \sin \sigma_{j} x  \end{array} \right)  , ~~ j\in\mathbb{N}.
$$

The fact that the  eigenfunctions $\{\Psi_{-1j},\Psi_{1j}, \Psi_{-2j}, \Psi_{2j}\}_{ j\in \mathbb{N}}$ are mutually orthogonal and  form a basis of $\mathrm H$  follows from the fact that  $\mc A$ is skew-symmetric and has  a compact resolvent (Lemma \ref{skew-adjoint}).

%    \begin{eqnarray} \nonumber && \{z_{1j}=\frac{1}{\lambda_{2j}^+}\sin \sigma_{j} x,~~ z_{2j}=\frac{b_2}{\lambda_{2j}^+}\sin \sigma_{j} x, ~~ z_{3j}=\sin \sigma_{j} x,~~ z_{4j}=b_2   \sin \sigma_{j} x,~~ j\in\mathbb{N}\},\\
%\nonumber  && \{z_{1j}=\frac{1}{\lambda_{2j}^+}\sin \sigma_{j} x, ~~ z_{2j}=\frac{b_2}{\lambda_{2j}^+}\sin \sigma_{j} x, ~~ z_{3j}=-\sin \sigma_{j} x,~~ z_{4j}=-b_2   \sin \sigma_{j} x,~~ j\in\mathbb{N}\}.
% \end{eqnarray}

%Finally, we only prove (\ref{norm-init}) since (\ref{sol}) and (\ref{init-F}) follow from (\ref{EVs}) and (\ref{evectors}).
Finally, prove (\ref{norm-init}).
Since \begin{eqnarray}\label{lazim}\rho+ b_1^2\mu=\zeta_1^2\left(\alpha_1 + \beta(\gamma-b_1^2)\right), \quad  \rho+ b_2^2\mu=\zeta_2^2\left(\alpha_1 + \beta(\gamma-b_2^2)\right),
\end{eqnarray}
%we obtain
%\begin{eqnarray}\nonumber \|\Phi_m(x,0)\|^2_{\mathrm{H}}&=&  L\left(2\mu + \rho \left(\frac{1}{b_1^2}+\frac{1}{b_2^2}\right)\right)\\
%\label{guzel}  &=& O(1)  \quad {\rm as}\quad   {\tilde q_m}\to \infty.
%\end{eqnarray}
by (\ref{inner}) and (\ref{ben}),  a direct calculation leads to
\begin{eqnarray}
\nonumber \|\Psi^0\|_{\mathrm{H}}^2 &=& \frac{L}{2} \left( \sum\limits_{j \in \mathbb{N}} (\rho + b_1^2 \mu + \alpha_1\zeta_1^2 + \beta (\gamma-b_1)^2\zeta_1^2 )\left(|c_{1j}|^2 + |d_{1j}|^2\right)\right.\\
\nonumber && \quad\quad\quad\quad  \left.+(\rho + b_2^2 \mu + \alpha_1\zeta_2^2 + \beta (\gamma-b_2)^2\zeta_2^2 )\left(|c_{2j}|^2 + |d_{2j}|^2\right)\right)\\
\nonumber &=&L\left( \sum\limits_{j \in \mathbb{N}} (\rho + b_1^2 \mu)\left(|c_{1j}|^2 + |d_{1j}|^2\right)+ (\rho + b_2^2 \mu )\left(|c_{2j}|^2 + |d_{2j}|^2\right)\right).\end{eqnarray}
Setting $\tilde C_1=L\min{\left(\rho + b_1^2 \mu, \rho + b_2^2 \mu \right)},$  $\tilde C_2=L\max{\left(\rho + b_1 ^2\mu ^2, \rho + b_2^2 \mu \right)}$ results in (\ref{norm-init}). $\square$

%%%%%%%%%%%%%%%%%%%

\end{document}